\documentclass[11pt]{article}

\RequirePackage[OT1]{fontenc}
\usepackage{fullpage}

\usepackage{graphicx}
\usepackage{lineno,hyperref}
\usepackage{amsmath,amssymb, amsbsy}
\usepackage{mathrsfs}
\usepackage{extarrows}
\usepackage{multirow,epsfig}
\usepackage{float}
\usepackage{xcolor}
\usepackage{booktabs}
\usepackage{graphicx}
\usepackage{subfigure}
\usepackage{caption}
\usepackage{epstopdf}
\usepackage[authoryear]{natbib}

\newtheorem{remark}{Remark}
\newtheorem{theorem}{Theorem}
\newtheorem{proposition}{Proposition}

\newtheorem{corollary}{Corollary}

\newcommand{\trans}{^{\mbox{\tiny{T}}}}
\newcommand{\bSig}{\boldsymbol{\Sigma}}
\newcommand{\hbSig}{\widehat{\boldsymbol{\Sigma}}}
\newcommand{\bOme}{\boldsymbol{\Omega}}
\newcommand{\hbOme}{\widehat{\boldsymbol{\Omega}}}
\newcommand{\tbOme}{\tilde{\boldsymbol{\Omega}}}
\newcommand{\bbeta}{\boldsymbol{\beta}}
\newcommand{\hbbeta}{\widehat{\boldsymbol{\beta}}}
\newcommand{\tbbeta}{\tilde{\boldsymbol{\beta}}}
\newcommand{\bh}{\boldsymbol{h}}
\newcommand{\be}{\boldsymbol{e}}

\newcommand{\mP}{\mathbf{P}}
\newcommand{\mE}{\mathbf{E}}
\newcommand{\mX}{\mathbf{X}}

\def \A{{\bf A}}
\def \B{{\bf B}}
\def\a{{\bf a}}

\def\r{{\bf r}}

\newcommand{\mR}{\mathbb{R}}
\DeclareMathOperator*{\argmin}{arg\,min} 
\newcommand{\tr}{\mbox{tr}}

\author{
	Zeyu Wu
	\footnote{Shanghai Jiao Tong University, Shanghai, China, Email: yufeng168@sjtu.edu.cn} ~
	Cheng Wang
	\footnote{Shanghai Jiao Tong University, Shanghai, China, Email: chengwang@sjtu.edu.cn} ~
	Weidong Liu
	\footnote{Shanghai Jiao Tong University, Shanghai, China, Email: weidongl@sjtu.edu.cn}
}

\begin{document}
	\title{A unified precision matrix estimation framework via sparse column-wise inverse operator under weak sparsity}
	\author{
		Zeyu Wu
		\footnote{Shanghai Jiao Tong University, Shanghai, China, Email: yufeng168@sjtu.edu.cn} ~
		Cheng Wang
		\footnote{Shanghai Jiao Tong University, Shanghai, China, Email: chengwang@sjtu.edu.cn} ~
		Weidong Liu
		\footnote{Shanghai Jiao Tong University, Shanghai, China, Email: weidongl@sjtu.edu.cn}
	}
	
	\date{}
	\maketitle
	\begin{abstract}
		In this paper, we estimate the high dimensional precision matrix under the weak sparsity condition where many entries are nearly zero.  We revisit the sparse column-wise inverse operator (SCIO) estimator \cite{liu2015fast} and derive its general error bounds under the weak sparsity condition. A unified framework is established to deal with various cases including the heavy-tailed data, the non-paranormal data, and the matrix variate data. These new methods can achieve the same convergence rates as the existing methods and can be implemented efficiently.  
	\end{abstract}
	\noindent
	\textbf{keywords:}
	Gaussian graphical model, high dimensional data, Lasso, precision matrix, weak sparsity

\section{Introduction}

In the high dimensional data analysis, estimating the population covariance matrix $\bSig$ and the precision matrix $\bOme=\bSig^{-1}$ are fundamental problems. The covariance matrix or the precision matrix characterizes the structure of the relation among covariates. Many statistical references can benefit from these structures if they can be precisely estimated.  We refer to \cite{tong2014estimation}, \cite{fan2016overview} and \cite{cai2017global} for recent reviews.

In the high dimensional regime where the data dimension is very large, estimating the covariance matrix or the precision matrix is challenging since the freedom of parameters are squared order of the data dimension. In literature, there have been a variety of methods which are proposed to estimate $\bSig$ or $\bOme$. For the population covariance matrix $\bSig$, \cite{bickel2008covariance} and \cite{el2008operator} constructed consistent estimation by thresholding the sample covariance matrix. This idea was further developed by \cite{rothman2009generalized} and \cite{cai2011adaptive}. For the precision matrix estimation, there have been many methods based on penalization. \cite{yuan2007model} considered a $\ell_1$-penalized Gaussian maximum likelihood estimator and \cite{friedman2008sparse} designed an efficient block coordinate descent algorithm for this so-called graphical Lasso  method. Noting that estimating each column of the precision matrix can be formulated as a linear regression problem,  \cite{Cai2011A} proposed a constrained $\ell_1$ minimization estimator (CLIME) based on Dantzig selector \citep{candes2007dantzig} and a similar estimator was also introduced by \cite{yuan2010high}.    \cite{liu2015fast} further proposed a sparse column-wise inverse operator estimator which was essentially a Lasso-type analog of their earlier CLIME method.  \cite{zhang2014sparse} considered a symmetric loss function which yielded a symmetric estimator directly, whereas other methods  \citep[e.g.,][]{yuan2010high, Cai2011A, liu2015fast} needed an additional symmetrization step.   
There have been a huge number of papers addressing the problem of the covariance matrix and the precision matrix estimation. For related methods and their connections, see the book by \cite{wainwright2019high}(e.g., Chapters 6 and 11).  

Under high dimensional settings, to obtain a consistent estimator of the covariance matrix or the precision matrix,  a sparsity condition is often imposed on the true matrix. Namely, many entries of the matrix are exactly zero or nearly so.  
A target matrix  $\A=(a_{ij})\in \mR^{p\times p}$ is said to be (strong) sparse or $\ell_0$ sparse,  meaning that 
\begin{align*}
	\mbox{$\ell_0$~sparse:}~\max_{j=1,\ldots,p}   \sum_{i=1}^p  |a_{ij}|^0  \leq s.
\end{align*}
That is, each column has at most $s$ nonzero elements. The $\ell_0$ sparsity condition requires that most entries of the matrix are exactly zero. A natural relaxation is to consider the weak sparsity or $\ell_q$ sparsity, that is 
\begin{align*}
	\mbox{$\ell_q$~sparse:}~\max_{j=1,\ldots,p}   \sum_{i=1}^p  |a_{ij}|^q  \leq s_q,
\end{align*}
for some $q \in [0,1)$ and $s_q>0$ is a radius. Note that the $\ell_0$ sparsity condition is a special case of $\ell_q$ sparsity condition with $q=0$ and $s_q=s$. 
For estimating the covariance matrix $\bSig$, \cite{bickel2008covariance} firstly provided the consistency result under the $\ell_q$ sparsity condition. Later,  \cite{rothman2009generalized} and \cite{cai2011adaptive} also studied the estimator with the weak $\ell_q$ sparsity condition.  These estimators have explicit forms  based on thresholding the elements of the sample covariance matrix and the theoretical analysis is relatively straightforward. As far as the precision matrix, the theoretical analysis is more challenging since the estimator usually does not have an explicit form.  \cite{ravikumar2011high} firstly established convergence rates for the graphical Lasso. In details, under the strong sparse condition,  they derived the bounds under the matrix element-wise infinity, the spectral and the Frobenius norm.  The graphical Lasso,  SCIO \citep{liu2015fast},  together with the D-trace method \citep{zhang2014sparse} are all based on a loss function with a $\ell_1$ penalty term. Technically, they all used the primal-dual witness technique \citep{wainwright2009sharp} or its extensions to prove the consistency under the $\ell_0$ sparsity condition. Specially, in order to apply the primal-dual technique, an irrepresentability condition is necessary; see Assumption 1 of \cite{ravikumar2011high}, Section 5.2 of \cite{zhang2014sparse} and the formula (4) in \cite{liu2015fast}. From the perspective of variable selection or support recovery \citep{meinshausen2006high}, the strong sparsity condition is a reasonable and also important criteria to evaluate the estimation.  However, from the perspective of matrix estimation, this condition is too restricted. For example, the Toeplitz matrix $\bOme=(\rho^{|i-j|})_{p \times p}$ for some $\rho \in (-1,1)$ is an important matrix in statistics which does not satisfy the strong sparsity condition. 

In this work, we focus on the $\ell_q$ sparse or weak sparse case. Theoretically, for the Dantzig-type methods, \cite{yuan2010high}, \cite{Cai2011A} and \cite{cai2016estimating} studied the convergence bounds for the weak sparse matrices. Specially, they derived the error bounds under the spectral norm, the matrix element-wise infinity  norm and the Frobenius norm of the estimation. Under the weak sparse settings, there is few theoretical study on Lasso-type methods of the precision matrix estimation. Among them, \cite{sun2013sparse} and \cite{ren2015asymptotic} exploited the scaled lasso \citep{sun2012scaled} to establish the optimal convergence rate for their precision matrix estimators under the normality assumption.

In high dimensional data analysis, there are various specific situations where we need to estimate the precision matrix and many methods were proposed in literature. Usually, these methods are based on the well-known procedure CLIME or the graphical Lasso. For example, to study the data with heavy-tailed distributions, \cite{avella2018robust} considered a robust approach to estimate the population covariance matrix and also the precision matrix. In details, they used the robust covariance matrix as a pilot estimator and implemented the CLIME method. For non-Guassian data,  \cite{liu2009nonparanormal} introduced a non-paranormal graphical model to describe the correlation among covariates.   \cite{liu2012high} and \cite{xue2012regularized} proposed a non-parametric rank-based estimator to estimate the covariance matrix and the precision matrix of the non-paranormal graphical model.  In details,  \cite{liu2012high} proposed four precision matrix estimators by plugging the rank-based estimate into the Dantzig-type method \citep{yuan2010high}, the CLIME \citep{Cai2011A}, the graphical Lasso \citep{yuan2007model, friedman2008sparse} and the neighbourhood pursuit estimator\citep{meinshausen2006high}.  For the matrix valued data where the covariance matrix and the precision matrix have certain structures such as Kronecker product, \cite{leng2012sparse} and \cite{zhou2014gemini} proposed feasible methods to estimate the precision matrix for matrix data. In particular,  \cite{leng2012sparse} used  the graphical Lasso and \cite{zhou2014gemini} considered both the graphical Lasso and CLIME.  These methods are driven by CLIME or the graphical Lasso. Computationally, it is well known that the implementation of CLIME or the graphical Lasso is time-consuming. As our recent work \cite{wang2020efficient} showed,  for $p \gg n$, the computation complexities of SCIO and D-trace are $O(np^2)$ while the one of the graphical Lasso is $O(p^3)$ for general case \citep[e.g.,][Section 3]{witten2011new}. For CLIME, each column is a Dantzig-selector regression and the state of the art algorithm is ``flare" \citep{li2015flare} which is based on the linearized alternating direction method of multipliers proposed by \cite{wang2012linearized}. For each column, the subprogram involves a Lasso problem whose computation is time-consuming. A detailed comparison on the computation time of these methods can be found in Section 3. Motivated by the appealing computational efficiency, it is natural to ask whether we can establish comparable convergence rates for Lasso-type methods under the weak sparse case with mild conditions.

In this article, we revisit the SCIO method and generalize the theoretical properties of SCIO under the $\ell_q$ sparsity condition by a new analysis which is different from the proof of  \cite{liu2015fast}. In details, we exploit the oracle inequality of Lasso  \citep{ye2010rate, sun2012scaled} and get the basic inequality for the SCIO method. Therefore, we can derive error bounds from the basic inequality directly and relax the common irrepresentability condition which is necessary for the primal-dual witness technique. Accordingly, we provide a unified framework for the SCIO method based on different types of covariance matrix estimation under various cases. These new SCIO-based methods can get the same consistence results as the ones based on CLIME and can be implemented more efficiently.

The SCIO estimation can be regarded as an example of the general M-estimation  \citep[Chapter 9]{wainwright2019high}.   \cite{Negahban2012} provided a unified framework of M-estimators with decomposable regularizes including $\ell_1$ penalty. They obtained the error bounds under weak sparsity for several applications. However, the analysis of \cite{Negahban2012} is based on the restricted strong convexity (RSC) condition which holds for Gaussian or sub-Gaussian distributions. For other complicated cases such as the heavy-tailed or the non-paranormal assumption, it is challenging to verify the RSC condition.  In contrast, our analysis is based on the basic inequality and also the structure of the precision matrix such as the symmetrization procedure \citep{Cai2011A}. We can establish reliable convergence rates of the SCIO method under weak sparsity due to its neat structure and provide a unified framework which is applicable for various distribution cases including the heavy-tailed distribution and non-paranormal distribution.

The rest of the paper is organized as follows. In Section 2, we revisit the SCIO method and derive the non-asymptotic results under the $\ell_q$ sparsity condition. In Section 3, we consider four covariance matrix estimators: the common sample covariance matrix for sub-Gaussian data, a Huber-type estimator for heavy-tailed data \citep{avella2018robust}, the nonparametric correlation estimator for non-paranormal data \citep{liu2009nonparanormal}  and the estimator for matrix variate data \citep{zhou2014gemini}.  By plugging these estimators into the SCIO procedure, we can estimate the corresponding precision matrix and derive the error bounds under the weak sparsity condition. Moreover, we conduct simulations to illustrate the performance of these estimators. Finally, we provide some brief comments in Section 4 and all technical proofs of our theorems, propositions and corollaries are relegated to Appendix.

\section{Main results}
We begin with some basic notations and definitions. For a vector $\a=(a_1,\ldots,a_p)\trans \in \mR^p$, the vector norms are defined as follows
\begin{align*}
	|\a|_\infty=\max_{1\leq i \leq p}{|a_i|}, ~|\a|_1=\sum_{i=1}^{p}|a_i|, \mbox{and}~ |\a|_2=\sqrt{\sum_{i=1}^{p}a_i^2}.
\end{align*}
For a matrix $\A=(a_{ij})\in \mR^{p\times q}$, we define the matrix norms:
\begin{itemize}
	\item the element-wise $l_{\infty}$ norm $\|\A\|_{\infty}=\max_{1\leq i \leq p, 1 \leq j \leq q} |a_{ij}|$;
	\item the spectral norm $\|\A\|_2=\sup_{|\mathbf{x}|_2\leq 1}|\A \mathbf{x}|_2$;
	\item the matrix $\ell_1$ norm $\|\A\|_{L_1}=\max_{1\leq j\leq q}\sum_{i=1}^{p}|a_{ij}|$;
	\item the Frobenius norm $\|\A\|_F=\sqrt{\sum_{i=1}^p\sum_{j=1}^{q}a_{ij}^2}$;
	\item the element-wise $\ell_1$ norm $\|\A\|_1=\sum_{i=1}^p\sum_{j=1}^{q}|a_{ij}|$.
\end{itemize}
For index sets $J\subseteq \{1,\ldots, p\}$ and $K \subseteq \{1,\ldots, q\}$, $\A_{J,\cdot}$ and $\A_{\cdot,K}$ denote the sub-matrix of $\A$ with rows or columns whose indexes belong to $J$ or $K$ respectively. In particular, $\A_{i,\cdot}$ and $\A_{\cdot,j}$ are the $i$th row and $j$th column respectively.  For a set $J$, $|J|$ denotes the cardinality of $J$. For two real sequences $\{a_n\}$ and $\{b_n\}$,  write $a_n=O(b_n)$ if there exists a constant $C$ such that $|a_n|\leq C|b_n|$ holds for large $n$ and $a_n=o(b_n)$ if $\lim_{n \to \infty}a_n/b_n=0$. The constants $C,C_0,C_1,...$ may represent different values at each appearance.

\subsection{SCIO revisited}
Suppose that $\hbSig$ be an arbitrary estimator of the population covariance matrix $\bSig$. Taking $\hbSig$ as the common sample covariance matrix, \cite{liu2015fast} proposed the SCIO method which estimated the precision matrix column-wisely.  Let $\be_i$ be the $i$th column of a $p\times p$ identity matrix and write $\bOme=(\bbeta_1,\ldots,\bbeta_p)$. For each $i=1,\ldots,p$,  we have 
\begin{align*}
	\bSig \bbeta_i=\be_i,
\end{align*}
and SCIO estimates the vector $\bbeta_i$ via a $\ell_1$ penalized form 
\begin{align}\label{eq:1}
	\hbbeta_i=\argmin_{\bbeta\in \mR^{p}}\left\{\frac{1}{2} \bbeta\trans \hbSig \bbeta-\be_i \trans \bbeta+\lambda|\bbeta|_{1}\right\},
\end{align}
where $\lambda \geq 0$ is a tuning parameter. By stacking the resulting  $\hbbeta_i$ together, we can obtain the precision matrix estimator $\hbOme=(\hbbeta_1,\ldots,\hbbeta_p)$. Noting that $\hbOme$ may be asymmetric, a further symmetrization step is necessary. The final SCIO estimator $\tbOme=(\tilde{\omega}_{i j})_{p \times p}$ is defined as
\begin{align}
	\tilde{\omega}_{i j}=\tilde{\omega}_{j i}=\left\{ \begin{array}{ll}
		\hat{\beta}_{i j} & \mbox{if $|\hat{\beta}_{i j}|\leq |\hat{\beta}_{j i}|$};\\
		\hat{\beta}_{ji} & \mbox{otherwise.}\end{array} \right.
	\label{eq:2}
\end{align}

From the method's perspective, SCIO is closely related to other popular precision matrix estimation methods. For example, the dual problem of \eqref{eq:1} is a Dantzig-type optimization problem
\begin{align}
	\min |\bbeta|_{1} \quad \text { subject to } \quad|\hbSig \bbeta-\be_i|_{\infty} \leq \lambda,
	\label{eq:3}
\end{align}
which is exactly the CLIME method \citep{Cai2011A}.  In a matrix form, \eqref{eq:1} is equivalent to
\begin{align*}
	\hbOme={\argmin_{\B \in \mR^{p \times p}}}\left\{\frac{1}{2}\tr( \B \hbSig \B \trans)-\tr(\B) +\lambda \|\B\|_1 \right\}.
\end{align*}
To deal with the problem that the objective function above is not symmetric about $\B$, \cite{zhang2014sparse} proposed the D-trace method which used the loss function 
\begin{align*}
	\frac{1}{4} \tr(\B \hbSig \B \trans)+\frac{1}{4} \tr(\B \trans \hbSig \B) -\tr(\B).
\end{align*}
Computationally, SCIO and D-trace use quadratic loss functions which can be solved efficiently via standard optimization algorithms.  More details can be found in our recent work \citep{wang2020efficient}.

From the theoretical perspective,  SCIO \citep{liu2015fast},  together with the graphical Lasso  \citep{ravikumar2011high} and D-trace\citep{zhang2014sparse} are all based on a loss function combined with a $\ell_1$ penalty term. To show the consistency of the estimation, they all used the primal-dual witness technique \citep{wainwright2009sharp} or its extensions under the $\ell_0$ sparsity condition. To apply the primal-dual technique, an irrepresentability condition is also necessary; see Assumption 1 of \cite{ravikumar2011high}, Section 5.2 of \cite{zhang2014sparse} and the formula (4) in \cite{liu2015fast}. In this work, we focus on exploring the theoretical property of the SCIO under the weak sparsity condition and relaxing the irrepresentability condition.

\subsection{SCIO for weak sparsity}
Note that \cite{ye2010rate} and \cite{raskutti2011minimax} considered the Lasso  under the $\ell_q$ sparsity condition. Here we study SCIO under the weak $\ell_q$ sparsity condition.  

Before presenting the error bounds, we define the $s_p$-sparse matrices class: 
\begin{align*}
	\mathcal{U}_q(s_p,M_p)=\left\{\bOme=(\omega_{ij})_{p\times p}:\bOme\succ0 , \max _{1 \leq j \leq p} \sum_{i=1}^{p} {\left|\omega_{ij}\right|^q} \leq s_{p}, \|\bOme\|_{L_{1}} \leq M_{p}\right\},
\end{align*}
where $s_p$ represents the sparsity for columns of the precision matrix, and $M_p$ may grow with the data dimension $p$. This matrix class was defined by \cite{Cai2011A} and see also \cite{bickel2008covariance} for a similar definition about the population covariance matrix. 

The first theorem refers to a convergence bound for the optimization problem \eqref{eq:1} in vector norms. This bound is established under the $s_p$-sparse matrices class and stated from a non-asymptotic viewpoint.
\begin{theorem}\label{thm:1}
	Suppose $\bOme\in\mathcal{U}_q(s_p,M_p)$ for some $0 \leq q <1$.  Assume that $\lambda\geq 3\|\bOme\|_{L_1}\|\hbSig-\bSig\|_\infty$ and $\|\bOme\|_{L_1}^{-q}\lambda^{1-q}s_p\leq \frac{1}{2}$ hold. Then 
	\begin{align*}
		&|\hbbeta-\bbeta^*|_1\leq 16\|\bOme\|^{1-q}_{L_1}s_p\lambda^{1-q},\\
		&|\hbbeta-\bbeta^*|_\infty\leq 4\|\bOme\|_{L_1}\lambda.
	\end{align*}
\end{theorem}
The technical proof of Theorem \ref{thm:1} is based on a basic inequality analogous to the one proposed by \cite{sun2012scaled}. We actually only take the special case where $w=\bbeta^*$ in \cite{sun2012scaled} into consideration.  \cite{ren2015asymptotic} used a similar trick and defined another weak sparsity based on a capped-$\ell_1$ measure. Their analysis depends on the Gaussian assumption and each element estimation $\hat{\omega}_{i j}$ needs to solve a scaled Lasso problem while SCIO can be solved  efficiently and yields the precision matrix estimation directly \citep{wang2020efficient}.

\begin{remark}
	In \cite{sun2013sparse} and \cite{ren2015asymptotic}, they used an alternative definition of the weak sparsity, i.e.,  the capped $\ell_1$ measure which is defined as $s_{t}(\bOme)=\max_{j}\sum_{i=1}^{p}\min\{1,|\beta^{*}_{ij}|/t\}$ for a threshold parameter $t$. See also the expository paper by \cite{Cai2016ejs}.  For every column $j$, by taking the index set $J=\{j~|~|\beta_{ij}|>t\}$, we have 
	\begin{align*}
		\sum_{i=1}^{p}\min\{1,|\beta^*_{ij}|/t\}=|J|+|\bbeta^{*}_{J^c}|/t.
	\end{align*}
	A slight modification of the proof in Theorem \ref{thm:1} yields
	\begin{align*}
		|\hbbeta-\bbeta^*|_1\leq\max\{12|\bbeta^*_{J^c}|_1, 16\lambda|J|\|\bOme\|_{L_1}\}\leq 16ts_{t}(\bOme)
	\end{align*}
	with $t=\lambda\|\bOme\|_{L_1}$. Thus we can extend our theoretical result to the capped $\ell_1$ measure. Moreover, by the discussion in \cite{ren2015asymptotic}, our result can also be extended to the weak $\ell_q$ ball sparsity condition.
\end{remark}

Given the non-asymptotic bounds of $\hbbeta-\bbeta^*$,  we remark that the symmetrization step \eqref{eq:2} is also crucial for the precision matrix estimation. \cite{yuan2010high} conducted another symmetrization procedure which was based on an optimization problem.   
In the next theorem, we present a non-asymptotic bound between the symmetric SCIO estimator $\tbOme$ and the true precision matrix $\bOme$ under the matrix norms.
\begin{theorem}\label{thm:2}
	Suppose $\bOme\in\mathcal{U}_q(s_p,M_p)$ for some $0 \leq q <1$.  Assume that $\lambda\geq 3\|\bOme\|_{L_1}\|\hbSig-\bSig\|_\infty$ and $\|\bOme\|_{L_1}^{-q}\lambda^{1-q}s_p\leq \frac{1}{2}$ hold. Then 
	\begin{align} \label{eq:4}
		\|\tbOme-\bOme\|_\infty\leq4\|\bOme\|_{L_1}\lambda,
	\end{align}
	and
	\begin{align}\label{eq:5}
		&\|\tbOme-\bOme\|_{L_1}\leq 66 ( \lambda \|\bOme\|_{L_1} )^{1-q}s_p.
	\end{align}
\end{theorem}

In Theorem \ref{thm:2}, we develop a unified framework for establishing convergence rates for the SCIO method. For any covariance matrix estimator $\hbSig$ with the bound $\|\hbSig-\bSig\|_\infty=o_p(1)$, the error bounds for the SCIO estimator under the matrix $\ell_\infty$ norm and the matrix $\ell_1$ norm are provided. Some specific examples with different choices of $\hbSig$ will be discussed in the next section. It is noted that we can further refine the error bound in Theorem \ref{thm:2}  by considering the tuning parameter $\lambda \geq 3\|\hbSig\bOme-\mathbf{I}\|_{\infty}$. However, for some covariance matrix estimators $\hbSig$, it is not trivial to derive the bound $\|\hbSig\bOme-\mathbf{I}\|_{\infty}$. For the conciseness and uniformity of our result statement, we consider $\lambda\geq 3\|\bOme\|_{L_1}\|\hbSig-\bSig\|_\infty$ here and will provide some comments about achieving the optimal bound for detailed applications later. Moreover, the $\ell_\infty$ norm \eqref{eq:4} and the matrix $\ell_1$ norm  \eqref{eq:5} are very useful and can yield other matrix bounds directly. For example, by the Gershgorin circle theorem, we can obtain the bound for the spectral norm
\begin{align*}
	\|\tbOme-\bOme\|_{2} \leq \|\tbOme-\bOme\|_{L_1} \leq 66 ( \lambda \|\bOme\|_{L_1} )^{1-q}s_p,
\end{align*}	
and for the Frobenius norm, we have
\begin{align*}
	\frac{1}{p}\|\tbOme-\bOme\|_{F}^{2} &\leq \|\tbOme-\bOme\|_{L_1}\|\tbOme-\bOme\|_{\infty}\\
	&\leq 264 \lambda^{2-q}  \|\bOme\|^{1-q}_{L_1} s_p.
\end{align*}
The matrix $\ell_1$ norm \eqref{eq:5} also plays an important role in many statistical inference. For example, an appropriate matrix $\ell_1$ bound can help to establish the asymptotic distribution of the test statistics in  \cite{cai2014two} or lead to the consistency of the thresholding estimation in \cite{wang2019cumulative}.

\section{Applications of the unified framework}
To illustrate the non-asymptotic bounds of Theorem \ref{thm:2}, we apply the SCIO method to several covariance matrix estimations. In details, for each plug-in covariance estimator $\hbSig$, we derive the bound $\|\hbSig-\bSig\|_\infty$ under high probability and apply Theorem \ref{thm:2} to show the consistency of the final precision matrix estimators $\tbOme$.

\subsection{Sample covariance matrix}
As a motivating application, we study the sample covariance matrix
\begin{align*}
	\hbSig_1=\frac{1}{n} \sum_{k=1}^{n}\left(\mX_{k}-\bar{\mX}\right)\left(\mX_{k}-\bar{\mX}\right)\trans
\end{align*}
where $\mX_k \in \mR^p$, $k=1,\ldots,n$ are independent and identically distributed (i.i.d.) samples and $\bar{\mX}=n^{-1} \sum_{k=1}^{n} \mX_{k}$ is the sample mean. \cite{liu2015fast} analyzed the sample covariance matrix and derived the consistency of SCIO under the $\ell_0$ sparsity condition and also a related irrepresentable condition. With the aim at more general $\ell_q$ sparsity setting, we state the technical conditions as the following:
\begin{itemize}
	\item[(A1).] (\textbf{Sparsity restriction}) Suppose that  $\bOme\in\mathcal{U}_q(s_p,M_p)$ for a given $q\in[0,1)$, where $s_p$ and $M_p$ satisfy the following assumption:
	\begin{align*}
		s_{p}M_p^{1-2q}=o\left(\frac{n}{\log p}\right)^{\frac{1}{2}-\frac{q}{2}}.
	\end{align*}
	\item[(A2).] (\textbf{Exponential-type tails}) 
	\noindent Suppose that $\log p = o(n)$. There exist positive numbers $\eta>0$ and $K>0$ such that
	\begin{align*}
		\mE \exp (\eta\left(X_{i}-\mu_{i}\right)^{2}) \leq K,  
	\end{align*}
	$\text { for all } 1 \leq i \leq p$.
	\item[(A3).]  (\textbf{Polynomial-type tails})
	\noindent Suppose that $p\leq cn^{\gamma}$ for some $\gamma, c>0$ and
	\begin{align*}
		\mE|X_i-\mu_i|^{4\gamma+4+\delta}\leq K,
	\end{align*} 
	for some $\delta>0$ and all $1 \leq i \leq p$.
\end{itemize}
The condition (A1) is an analogue of the formula (3) in \cite{liu2015fast} which is for the special case $q=0$.  The  conditions (A2) and (A3) are regular conditions which are used to control the tail probability of the variables. See also the assumptions of \cite{Cai2011A}.  Under the conditions (A2) and (A3), \cite{liu2015fast} proved the following proposition:
\begin{proposition}[Lemma 1, \citealt{liu2015fast}]
	For a given $\tau>0$ and a sufficiently large constant $C$, we have 
	\begin{align*}
		\mP\left(\|\hbSig_1-\bSig\|_{\infty}\geq C\sqrt{\frac{\log{p}}{n}}\right)\leq O(p^{-\tau}),
	\end{align*}
	under the assumption (A2) or
	\begin{align*}
		\mP\left(\|\hbSig_1-\bSig\|_{\infty}\geq C\sqrt{\frac{\log{p}}{n}}\right)\leq O(p^{-\tau}+n^{-\frac{\delta}{8}}),
	\end{align*}
	under the assumption (A3).
	\label{prop:1}
\end{proposition}

With these results, for the estimator $\tilde{\bOme}_1$ obtained by plugging $\hbSig_1$ into SCIO, we are ready to state our main results under the $\ell_q$ sparsity setting.
\begin{corollary}
	Let $\lambda=C_0  \sqrt{\log p / n}$ with $C_0$ being a sufficiently large number. For $\bOme\in\mathcal{U}_q(s_p,M_p)$, under assumptions (A1) and (A2) or (A3), we have
	\begin{align*}
		&\|\tbOme_1-\bOme\|_{\infty} \leq C_1 M_{p}^2 \sqrt{\frac{\log p}{n}},\\
		&\|\tbOme_1-\bOme\|_{L_1} \leq C_2 s_p M_{p}^{2-2q}\left(\frac{\log p}{n}\right)^{\frac{1}{2}(1-q)},
	\end{align*}
	with probability greater than $1-O(p^{\tau})$ or $1-O(p^{-\tau}+n^{-\frac{\delta}{8}})$. Here $C_1,C_2$ are sufficiently large constants which only depend on $q, s_p, M_p, C_0, \eta,K, \delta$.
	\label{cor:1}
\end{corollary}

By plugging the sample covariance matrix into SCIO, the estimation \eqref{eq:1} is similar to the classical Lasso regression problem and the error bounds considered here are analogous to prediction error bounds of Lasso regression problem \citep[Theorem 7.20]{wainwright2019high}. We adopt a different analysis from the primal-dual witness technique considered in \cite{liu2015fast} and remove the irrepresentability condition to obtain the error bounds under the $\ell_q$ sparsity setting. Correspondingly, there is no variable selection consistency results since the notion of variable selection is ambiguous for the $\ell_q$ sparsity. Moreover, we have the following remarks.
\begin{remark}
	Comparing to the CLIME method \citep{Cai2011A}, we derive the same convergence rates under the same conditions. This verifies the dual relation between Lasso and the Dantzig selector. \cite{2009Simultaneous} showed this point for the regression model and the results here demonstrate that the Lasso-type method and the Dantzig-type method for the precision matrix estimation also exhibit similar behaviors. 
\end{remark}

\begin{remark}
	If we impose stronger conditions on the tail distribution of $\bOme X_i$, i.e., conditions (C2) and (C2*) in \cite{liu2015fast},  we can get
	\begin{align*}
		\mP\left(\max _{1 \leq i \leq p}|\hbSig_1 \bbeta_{i}^{*}-\be_{i}|_{\infty} \geq C \sqrt{\frac{\log p}{n}}\right) = O(p^{-\tau}),
	\end{align*}
	or
	\begin{align*}
		\mP\left(\max _{1 \leq i \leq p}|\hbSig_1 \bbeta_{i}^{*}-\be_{i}|_{\infty} \geq C \sqrt{\frac{\log p}{n}}\right) = O(p^{-\tau}+n^{-\frac{\delta}{8}}),
	\end{align*}
	where $\bbeta_i^*$ is the $i$-th column of the true precision matrix $\bOme$. Then with some additional efforts, the error bounds
	\begin{align*}
		&\|\tbOme_1-\bOme\|_{\infty} \leq C_1 M_{p} \sqrt{\frac{\log p}{n}},\\
		&\|\tbOme_1-\bOme\|_{L_1} \leq C_2 s_p M_{p}^{1-q}\left(\frac{\log p}{n}\right)^{\frac{1}{2}(1-q)},
	\end{align*}
	hold with probability greater than $1-O(p^{-\tau})$ or $1-O(p^{-\tau}+n^{-\frac{\delta}{8}})$. These convergence rates actually achieve the minimax rate for estimating the true precision matrix $\bOme\in\mathcal{U}_q(s_p,M_p)$. See \cite{cai2016estimating} for more details.
\end{remark}

The main motivation to study the SCIO for weak sparsity is that it is computationally more efficient than other methods such as CLIME or the graphical Lasso \citep[e.g., Table 1 of][]{wang2020efficient}. Here, we further conduct several simulations to compare the computation time of these methods. In details, we include the scaled Lasso method(SLasso) which is implemented with the R package ``scalreg" provided by \cite{sun2013sparse},  the CLIME method which is implemented with the R package ``flare" developed by  \cite{li2015flare}, the graphical Lasso (gLasso) which is implemented with the R packages ``gLasso", ``BigQuic" or ADMM algorithm \citep[Section 6.5]{boyd2011distributed}, the D-trace and the SCIO which are implemented with the R package  ``EQUAL" developed by \cite{wang2020efficient}. The SLasso uses the default setting and for all other methods, the computation time is recorded in seconds and averaged over 5 replications on a solution path with 50 $\lambda$ values ranging from $\lambda_{max}$ to $\lambda_{max} \sqrt{\log{p}/n} $. Here $\lambda_{max}$ is the maximum absolute off-diagonal elements of the sample covariance matrix. All methods are evaluated on an Intel Core i7 3.3GHz and under R version 4.2.1 with an optimized BLAS implementation for Mac hardware. Table \ref{tab1} summarizes the computation time. Although the stopping criteria is different for each method, we can see from Table \ref{tab1} the superior efficiency of the SCIO method. 

\begin{table}[!htbp] 
	\caption{The average computation time (standard deviation) of the precision matrix estimation with $n=200$. } 
	\label{tab1} 
	\resizebox{\textwidth}{!}{%
		\begin{tabular}{@{\extracolsep{5pt}} ccccccc} 
			\\[-1.8ex]\hline 
			\hline \\[-1.8ex] 
			& $p$=50  	& $p$=100 & $p$=200 & $p$=400 & $p$=800 & $p$=1600 \\ 
			\hline \\[-1.8ex] 
			&\multicolumn{6}{c}{Case 1:  $\bOme=(0.5^{|i-j|})_{p \times p}$}\\  \\[-1.8ex] 
			SLasso(scalreg) & 0.91(0.08) & 4.96(0.14) & 95.56(0.45) & 191.35(1.96) &  0.00(0.00) &   0.00(0.00) \\ CLIME(flare) & 0.67(0.22) & 3.01(0.08) & 65.58(0.70) & 186.90(0.10) &  0.00(0.00) &   0.00(0.00) \\ gLasso(gLasso) & 0.04(0.00) & 0.25(0.00) &  1.78(0.02) &  10.60(0.40) & 64.81(0.41) & 737.06(2.21) \\ gLasso(BigQuic) & 0.36(0.03) & 0.76(0.02) &  2.19(0.03) &   6.25(0.01) & 21.14(0.06) &  92.73(0.36) \\ gLasso(ADMM) & 0.10(0.00) & 0.39(0.03) &  1.23(0.02) &   4.25(0.04) & 14.10(0.75) &  92.01(3.51) \\ D-trace(EQUAL) & 0.03(0.00) & 0.11(0.01) &  0.32(0.00) &   0.60(0.01) &  2.00(0.03) &   8.64(0.46) \\ SCIO(EQUAL) & 0.02(0.00) & 0.06(0.00) &  0.18(0.00) &   0.41(0.00) &  1.65(0.03) &   7.95(0.07) \\ 
			\hline \\[-1.8ex] 
			&\multicolumn{6}{c}{Case 2: $\bOme^{-1}=(0.5^{|i-j|})_{p \times p}$} \\  \\[-1.8ex] 
			SLasso(scalreg) & 0.84(0.15) & 5.42(0.26) & 117.17(3.41) & 183.69(7.32) &  0.00(0.00) &   0.00(0.00) \\ CLIME(flare) & 1.05(0.08) & 6.28(0.07) &  81.29(0.10) & 323.12(2.13) &  0.00(0.00) &   0.00(0.00) \\ gLasso(gLasso) & 0.05(0.00) & 0.24(0.00) &   1.48(0.05) &   9.66(0.28) & 79.82(2.82) & 762.80(7.78) \\ gLasso(BigQuic) & 0.40(0.01) & 1.00(0.05) &   2.71(0.02) &   8.71(0.07) & 39.17(0.54) & 212.43(5.25) \\ gLasso(ADMM) & 0.15(0.01) & 0.42(0.01) &   1.13(0.04) &   3.53(0.06) & 13.98(1.37) & 103.14(4.44) \\ D-trace(EQUAL) & 0.04(0.00) & 0.14(0.00) &   0.39(0.00) &   1.01(0.03) &  3.52(0.04) &  22.48(1.84) \\ SCIO(EQUAL) & 0.02(0.00) & 0.07(0.00) &   0.22(0.01) &   0.68(0.00) &  2.81(0.02) &  20.38(1.47) \\ 
			\hline \\[-1.8ex] 
	\end{tabular} }
\end{table} 

To further investigate the numerical performance of the SCIO estimation, we compare it with SLasso, gLasso, CLIME and D-trace. While the SLasso is tuning-free, we implement a five folds cross-validation procedure to select the tuning parameter $\lambda$ for all other methods.  The tuning parameter $\lambda$ is selected from 50 different values by minimizing the quadratic loss
\begin{align*}
	\operatorname{Loss}(\bOme)=\frac{1}{2}\tr( \bOme \hbSig \bOme \trans)-\tr(\bOme)
\end{align*}
where $\bOme$ is computed based on the training sample and $\hbSig$ is the sample covariance matrix of the test sample. To alleviate the bias of $\ell_1$ penalty, we also include the relaxed version \citep{meinshausen2006high, hastie2020best} of the estimator where a two-stage refitted estimator is obtained based on the support of the original estimator. These estimators are denoted by SLasso-R, gLasso-R, CLIME-R, D-trace-R and SCIO-R.  Table \ref{tab2} presents the estimation error for dimensions $p=100,200,400$ based on 100 replications. From the simulation results of Table \ref{tab1} and Table \ref{tab2}, we can conclude that SCIO enjoys comparable statistical convergence rates with superior computational efficiency in comparison to existing methods.

\begin{table}[!htbp] \centering 
	\caption{The average statistical error (standard deviation) of the precision matrix estimation with $n=200$.} 
	\label{tab2} 
	\resizebox{\textwidth}{!}{%
		\begin{tabular}{@{\extracolsep{5pt}} cccccccccc} 
			\\[-1.8ex]\hline 
			\hline \\[-1.8ex] 
			&\multicolumn{3}{c}{$p=100$}&\multicolumn{3}{c}{$p=200$}&\multicolumn{3}{c}{$p=400$}\\
			\hline \\[-1.8ex]
			& Spectral & Frobenius & $L_1$ & Spectral & Frobenius & $L_1$ & Spectral & Frobenius & $L_1$ \\ 
			\hline \\[-1.8ex] 
			&\multicolumn{9}{c}{Case 1:  $\bOme=(0.5^{|i-j|})_{p \times p}$}\\  \\[-1.8ex]
			SLasso &  1.60(0.03) &  5.39(0.07) &  1.97(0.06) &  1.68(0.02) &  8.06(0.07) &  2.04(0.05) &  1.75(0.02) & 11.91(0.07) &  2.12(0.04) \\ 
			SLasso-R &  1.20(0.05) &  4.42(0.11) &  2.14(0.16) &  1.30(0.04) &  6.76(0.12) &  2.50(0.21) &  1.37(0.02) & 10.08(0.09) &  2.75(0.22) \\ 
			gLasso &  1.76(0.03) &  5.99(0.12) &  2.43(0.06) &  1.91(0.02) &  9.34(0.08) &  2.52(0.08) &  2.00(0.01) & 13.94(0.08) &  2.59(0.07) \\ 
			gLasso-R &  1.49(0.06) &  4.86(0.09) &  1.94(0.13) &  1.54(0.06) &  6.98(0.11) &  2.06(0.13) &  1.60(0.05) & 10.03(0.14) &  2.23(0.12) \\ 
			CLIME &  1.85(0.03) &  6.22(0.10) &  2.19(0.04) &  1.82(0.02) &  8.58(0.09) &  2.24(0.07) &  1.96(0.02) & 13.49(0.10) &  2.32(0.04) \\ 
			CLIME-R &  1.40(0.03) &  4.60(0.09) &  1.80(0.13) &  1.51(0.06) &  6.78(0.13) &  2.06(0.11) &  1.29(0.10) &  6.94(0.23) &  1.99(0.16) \\ 
			D-trace &  1.61(0.05) &  5.30(0.13) &  2.04(0.06) &  1.75(0.03) &  8.29(0.14) &  2.13(0.04) &  1.86(0.02) & 12.60(0.18) &  2.20(0.03) \\ 
			D-trace-R &  1.46(0.06) &  4.80(0.09) &  1.89(0.15) &  1.51(0.04) &  6.87(0.09) &  2.00(0.10) &  1.59(0.06) &  9.87(0.11) &  2.17(0.13) \\ 
			SCIO &  1.61(0.04) &  5.35(0.13) &  2.03(0.07) &  1.79(0.02) &  8.55(0.12) &  2.13(0.04) &  1.92(0.02) & 13.12(0.09) &  2.24(0.04) \\ 
			SCIO-R &  1.44(0.05) &  4.77(0.08) &  1.84(0.12) &  1.47(0.04) &  6.83(0.08) &  1.92(0.11) &  1.52(0.04) &  9.76(0.09) &  2.06(0.11) \\ 
			
			&\multicolumn{9}{c}{Case 2: $\bOme^{-1}=(0.5^{|i-j|})_{p \times p}$} \\  \\[-1.8ex]
			SLasso &  0.78(0.11) &  3.37(0.15) &  1.10(0.15) &  0.83(0.07) &  5.03(0.11) &  1.15(0.10) &  0.90(0.07) &  7.51(0.14) &  1.24(0.11) \\ 
			SLasso-R &  1.20(0.15) &  4.23(0.23) &  2.18(0.34) &  1.34(0.13) &  6.46(0.23) &  2.48(0.33) &  1.49(0.10) &  9.46(0.27) &  2.84(0.34) \\ 
			gLasso &  1.03(0.05) &  4.72(0.12) &  1.56(0.09) &  1.10(0.05) &  7.02(0.13) &  1.82(0.10) &  1.18(0.04) & 10.80(0.17) &  2.04(0.12) \\ 
			gLasso-R &  0.83(0.17) &  2.56(0.20) &  1.08(0.22) &  0.91(0.11) &  3.66(0.20) &  1.20(0.19) &  1.02(0.16) &  5.24(0.28) &  1.36(0.25) \\ 
			CLIME &  1.02(0.07) &  4.34(0.17) &  1.39(0.09) &  1.11(0.06) &  6.63(0.18) &  1.51(0.08) &  1.27(0.05) & 11.86(0.14) &  1.65(0.07) \\ 
			CLIME-R &  0.80(0.16) &  2.43(0.17) &  1.03(0.21) &  0.87(0.10) &  3.48(0.18) &  1.12(0.15) &  0.97(0.12) &  4.92(0.24) &  1.24(0.17) \\ 
			D-trace &  0.91(0.07) &  3.84(0.21) &  1.26(0.10) &  1.00(0.05) &  5.95(0.20) &  1.37(0.08) &  1.09(0.04) &  9.05(0.17) &  1.46(0.06) \\ 
			D-trace-R &  0.80(0.15) &  2.44(0.18) &  1.01(0.17) &  0.84(0.14) &  3.46(0.18) &  1.08(0.19) &  0.96(0.14) &  4.90(0.24) &  1.22(0.20) \\ 
			SCIO &  0.89(0.07) &  3.84(0.20) &  1.21(0.09) &  1.02(0.04) &  6.16(0.20) &  1.34(0.07) &  1.14(0.04) &  9.76(0.17) &  1.47(0.06) \\ 
			SCIO-R &  0.79(0.14) &  2.45(0.18) &  1.02(0.18) &  0.85(0.16) &  3.51(0.19) &  1.17(0.27) &  0.97(0.16) &  4.99(0.22) &  1.32(0.23) \\ 
			\hline \\[-1.8ex] 
	\end{tabular}}
\end{table}

Next we conduct some simulations to illustrate the developed theoretical results. Firstly, in order to show that the irrepresentable condition is not necessary, we revisit the diamond graph example in \cite{ravikumar2011high} and consider a block precision matrix:
\begin{align*}
	\bOme=(\operatorname{diag}(\A,\cdots,\A))^{-1},
\end{align*}
where
\begin{align*}
	\A=
	\begin{pmatrix}
		1&\rho&\rho&2\rho^2 \\
		\rho&1&0&\rho\\
		\rho&0&1&\rho\\
		2\rho^2&\rho&\rho&1
	\end{pmatrix}
	\in\mR^{4 \times 4},
\end{align*}
and $\rho \in (-1/\sqrt{2},1/\sqrt{2})$ which ensures the positive definiteness of the covariance matrix.  The irrepresentable condition of the graphical Lasso  \citep{ravikumar2011high} holds for $|\rho|<(\sqrt{2}-1)/2$  and the irrepresentable conditions in \cite{liu2015fast} and \cite{meinshausen2006high} require that $|\rho|<1/2$.  

\begin{figure}
	\centerline{
		\begin{tabular}{ccc}				
			\psfig{figure=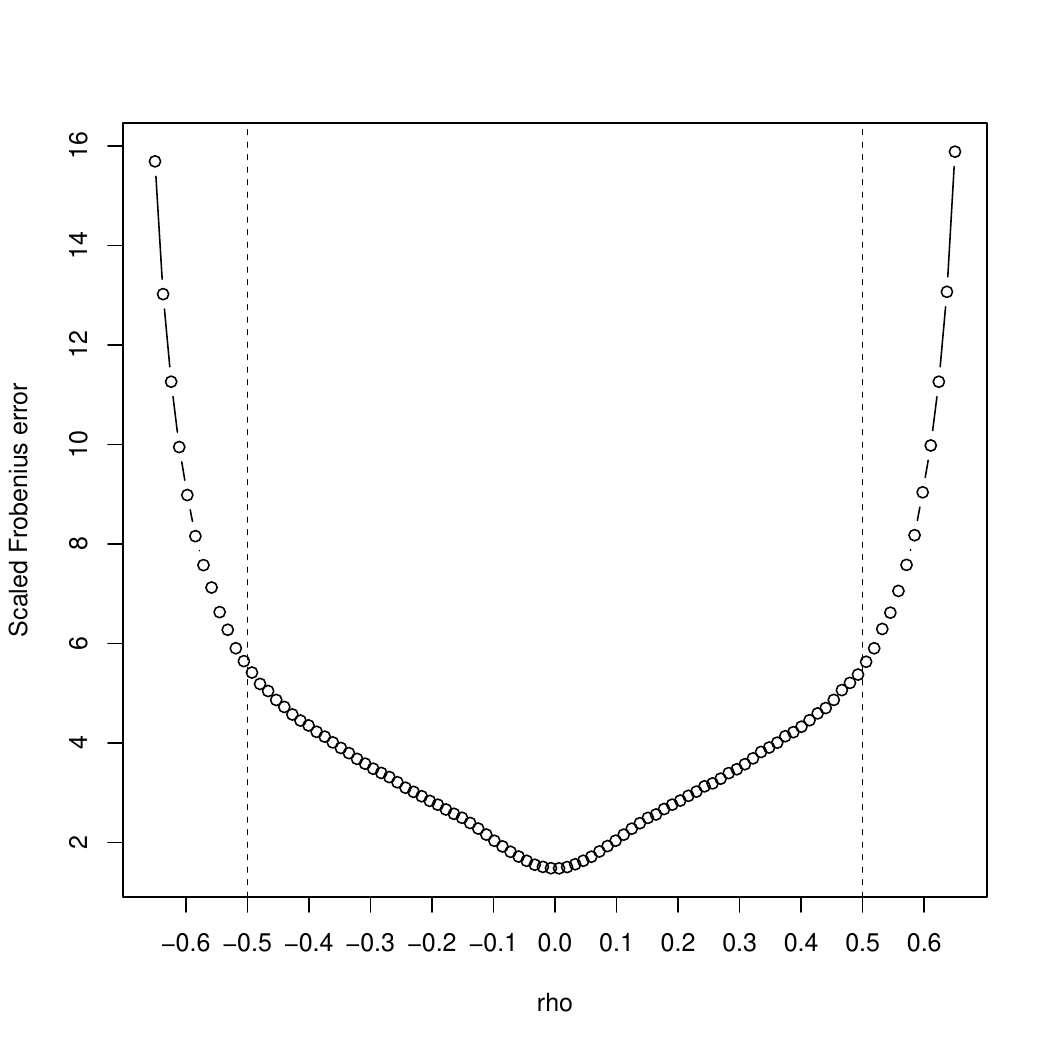, width=0.3\linewidth,angle=0} &
			\psfig{figure=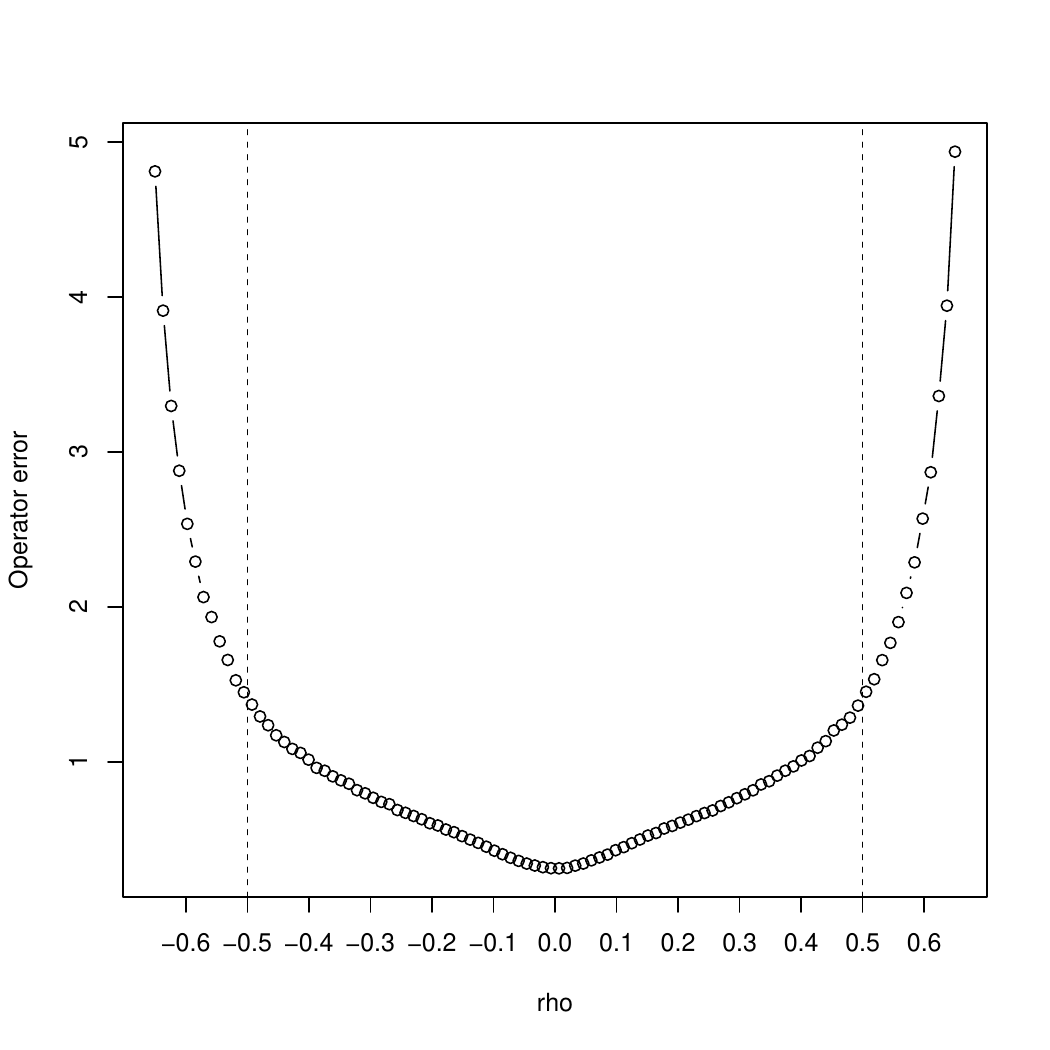, width=0.3\linewidth,angle=0}&
			\psfig{figure=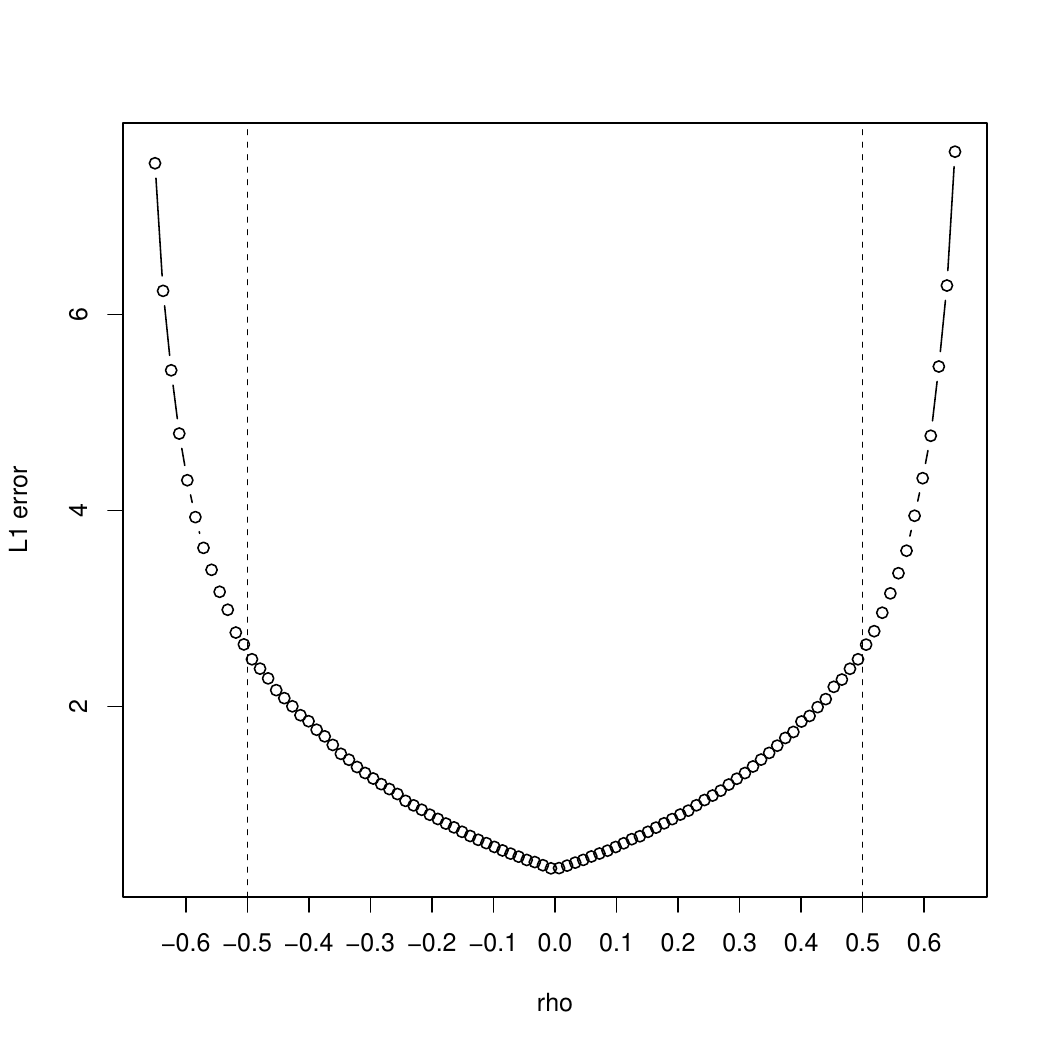, width=0.3\linewidth,angle=0}  \\
			(a) Frobenius norm & (b) Spectral norm & (c) $L_1$ norm  
		\end{tabular}
	}
	\caption{Plots of the estimation errors versus the parameter $\rho$ under three norms. The dash vertical lines indicate the boundaries of the irrepresentable condition. }
	\label{fig:1}
\end{figure}

Figure \ref{fig:1} shows the performance of the SCIO estimation for $\rho \in [-0.65,0.65]$. The sample is generated by the multivariate Gaussian distribution $\mathcal{N}_p(0,\bOme^{-1})$ where $p=100$ and $n=200$.  We plot the spectral norm, the matrix $\ell_1$ norm and the scaled Frobenius norm of $\bOme-\hbOme$. For the brevity, the tuning parameter $\lambda$ is chosen by minimizing the  matrix norms.  From these figures, we can observe that all the errors vary smoothly when $\rho$ is changing. Particularly, these errors do not drop drastically around the critical boundary value $|\rho|=0.5$. This phenomenon indicates that even though the validity of irrepresentable condition fails when $|\rho|\geq 0.5$, the performance of the SCIO estimation does not become worse drastically. Therefore, it is reasonable to relax the extra irrepresentable condition for SCIO.

To illustrate the consistent results for weak sparse cases,  we further conduct numerical studies where $\bOme=(\omega_{ij})_{p \times p}=(\rho^{|i-j|})_{p \times p}$ for some $\rho \in (0,1)$. For a fixed $q \in (0,1)$, we know
\begin{align*}
	\max _{1 \leq j \leq p} \sum_{i=1}^{p} {\left|\omega_{ij}\right|^q} \approx  1+2 \sum_{k=1}^\infty \rho^{kq}=\frac{1+\rho^q}{1-\rho^q}:=s_{p}.
\end{align*}
Hence the parameter $\rho$ measures the sparsity level of the true precision matrix. When $\rho$ is small, the decay phenomenon is salient and the matrix tends to be more sparse. When $\rho$ is large, the number of elements with small magnitude accounts for less proportion of all elements.  

\begin{figure}
	\centerline{
		\begin{tabular}{ccc}				
			\psfig{figure=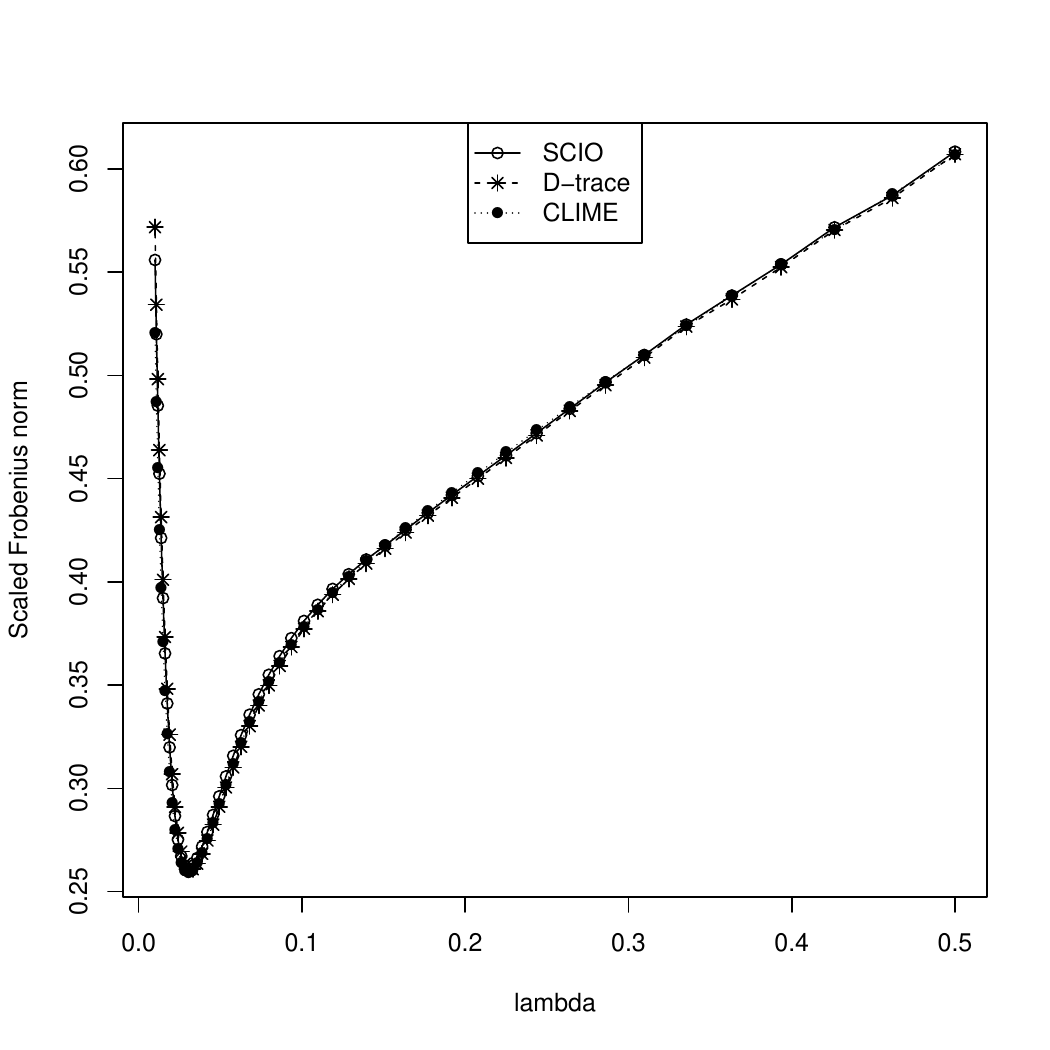, width=0.3\linewidth,angle=0} &
			\psfig{figure=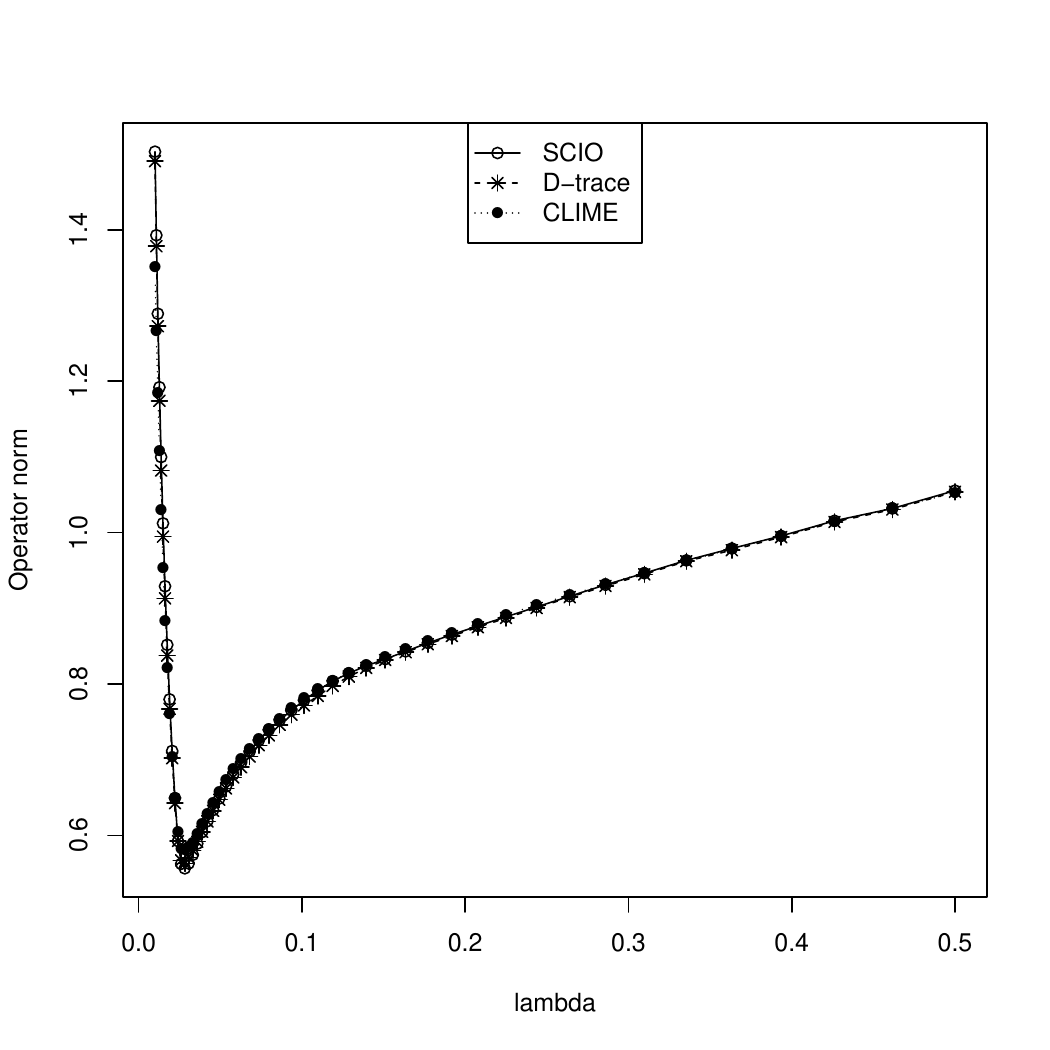, width=0.3\linewidth,angle=0}&
			\psfig{figure=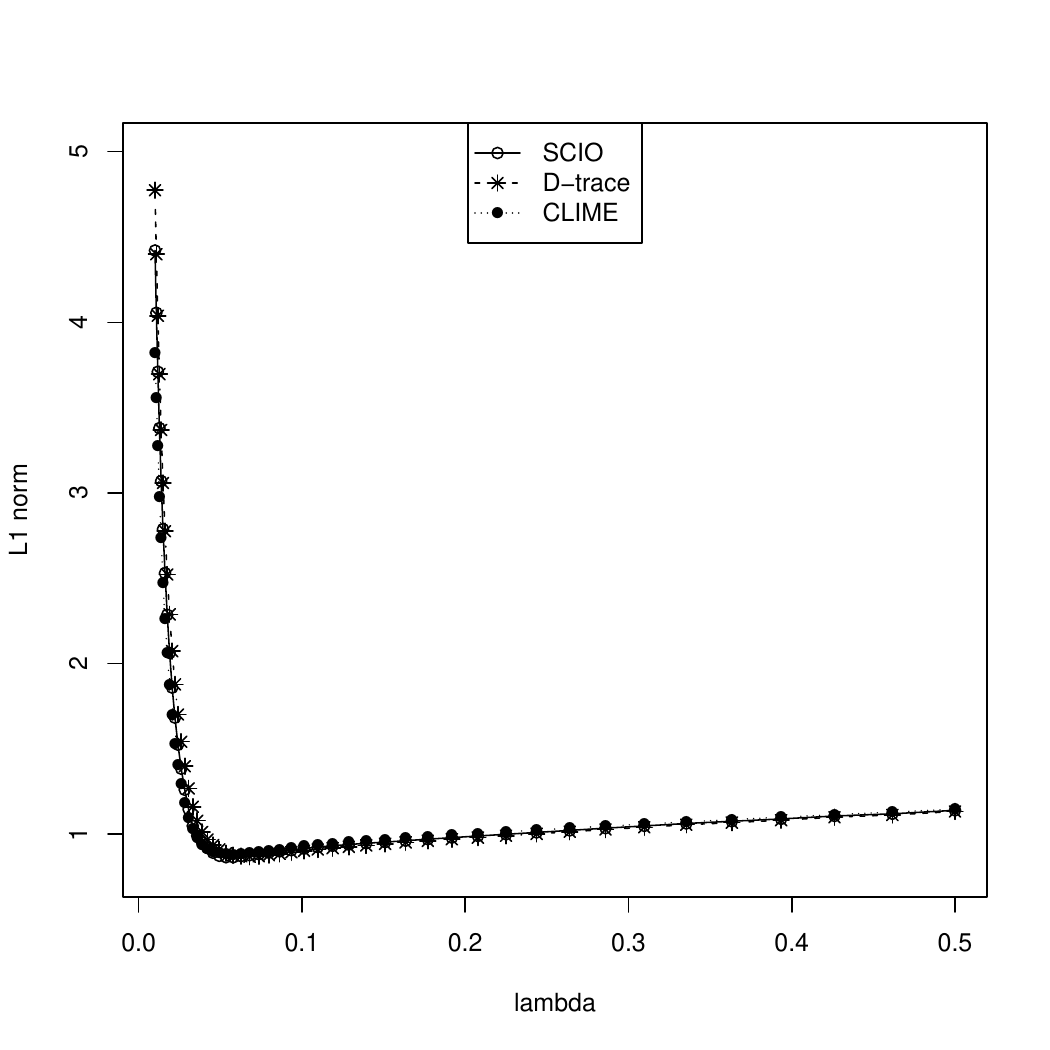, width=0.3\linewidth,angle=0}  \\
			(a) Frobenius norm for $\rho=0.2$ & (b) Spectral norm for $\rho=0.2$ & (c) $L_1$ norm for $\rho=0.2$\\
			\psfig{figure=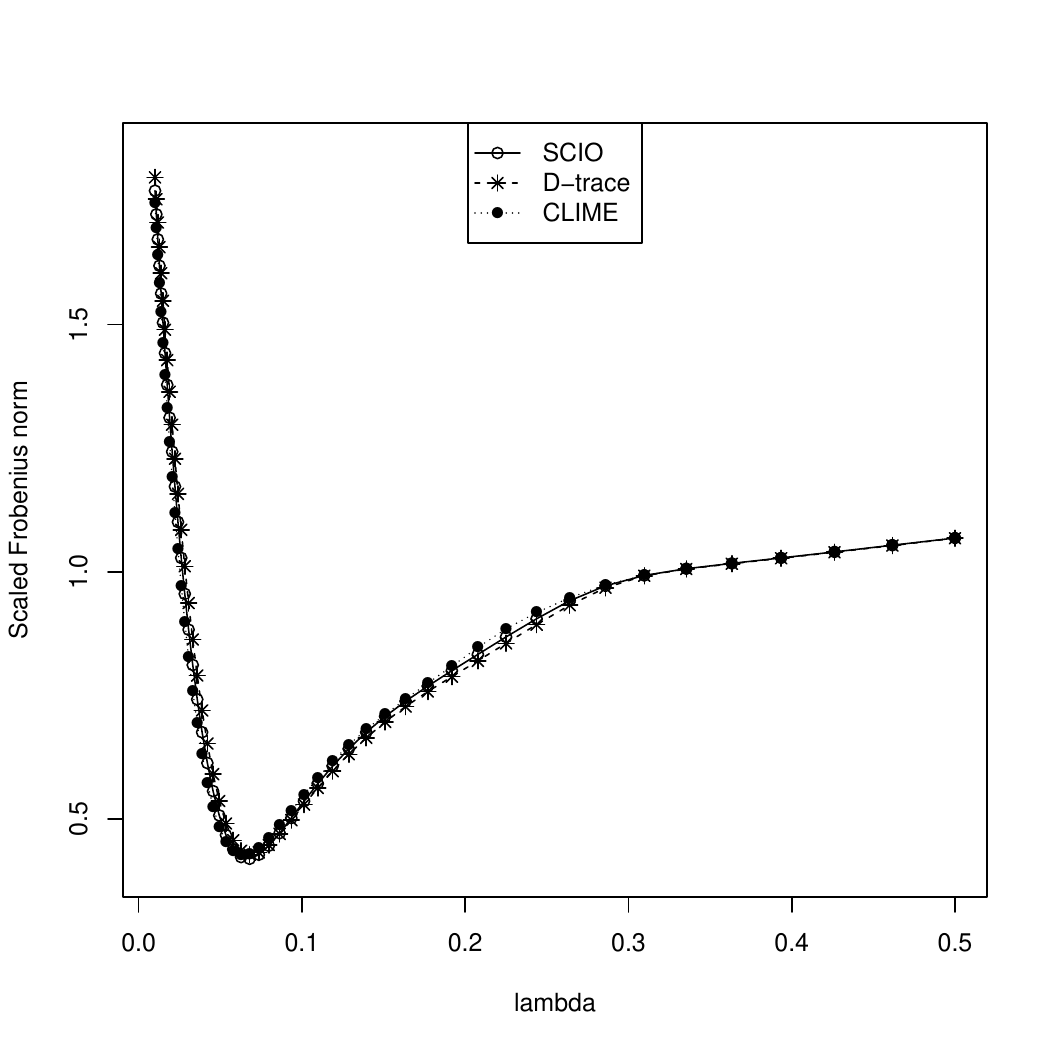, width=0.3\linewidth,angle=0} &
			\psfig{figure=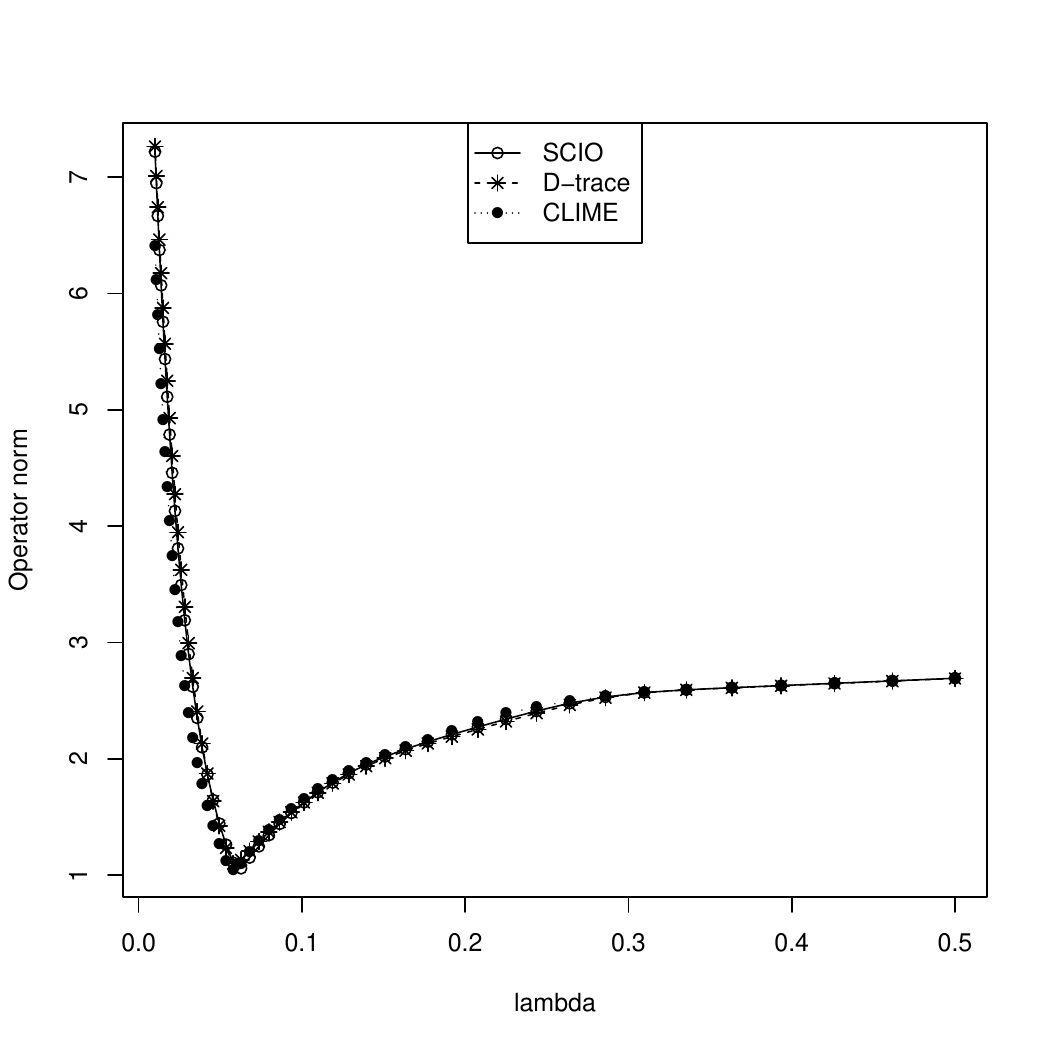, width=0.3\linewidth,angle=0}&
			\psfig{figure=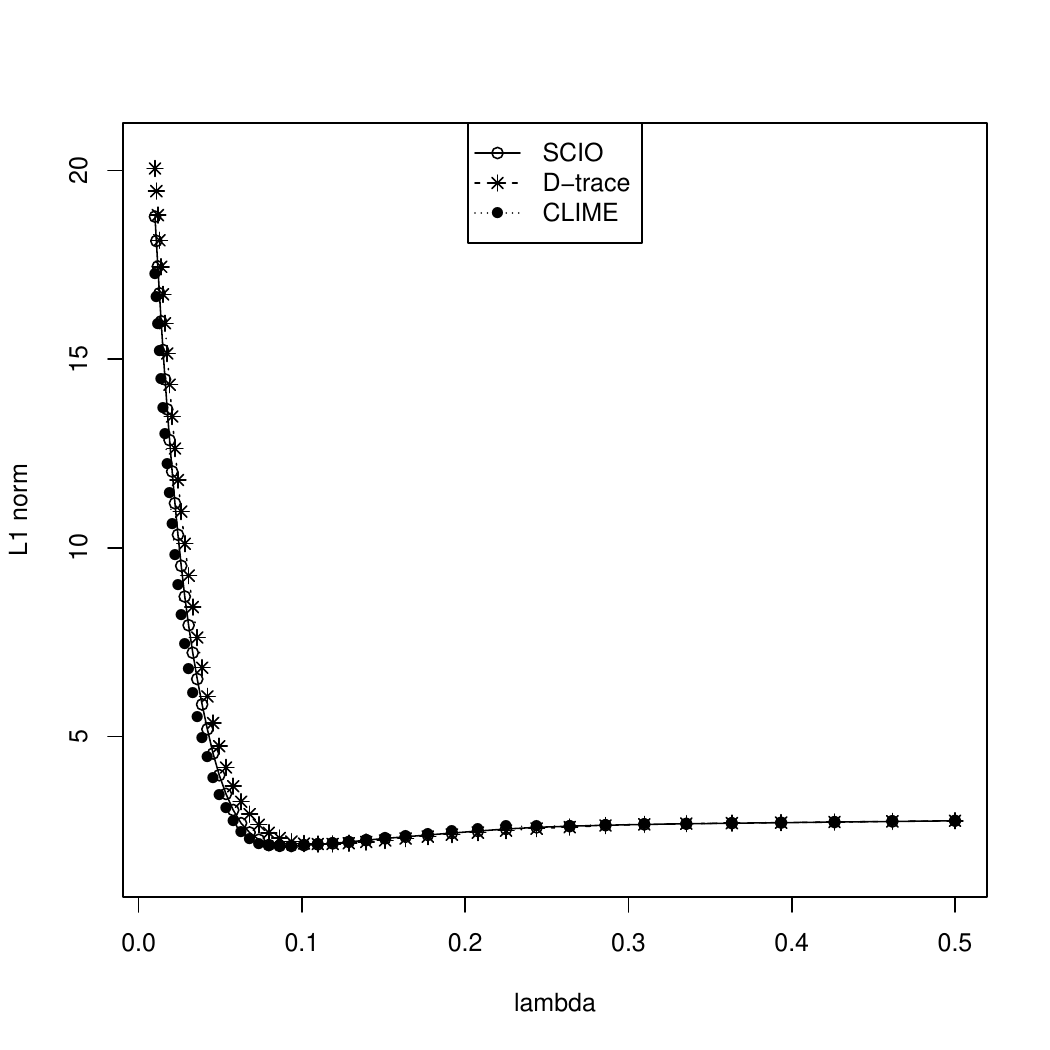, width=0.3\linewidth,angle=0}  \\
			(d) Frobenius norm for $\rho=0.5$ & (e) Spectral norm for $\rho=0.5$ & (f) $L_1$ norm for   $\rho=0.5$ \\
			\psfig{figure=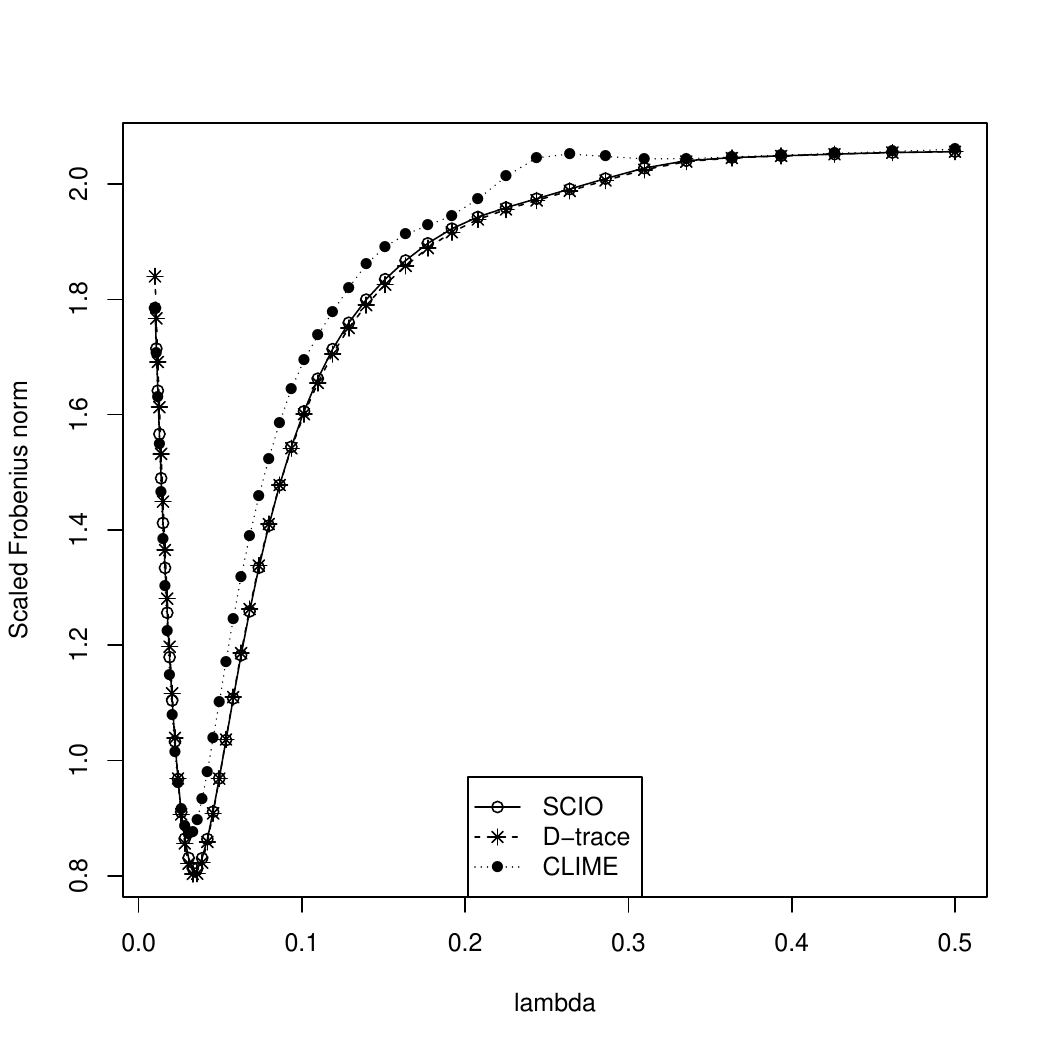, width=0.3\linewidth,angle=0} &
			\psfig{figure=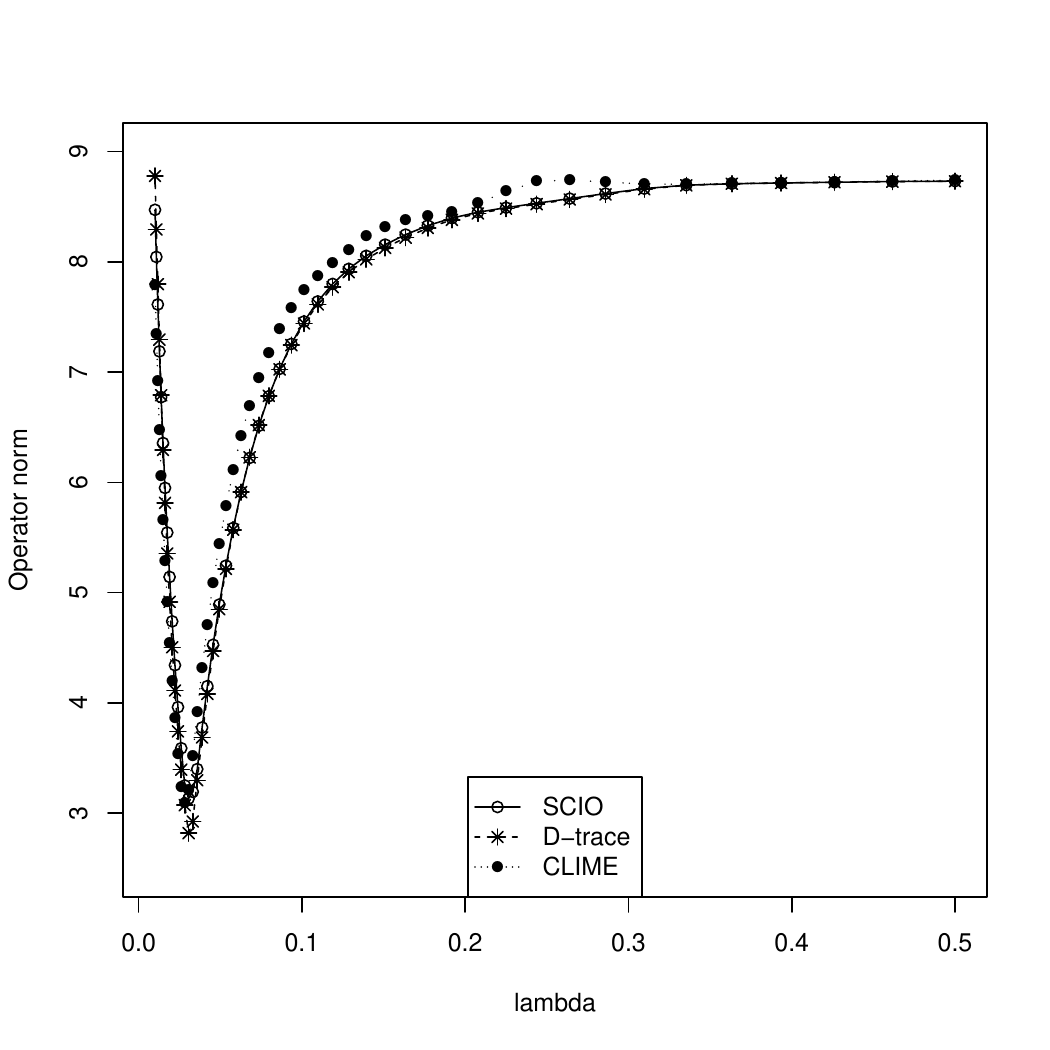, width=0.3\linewidth,angle=0}&
			\psfig{figure=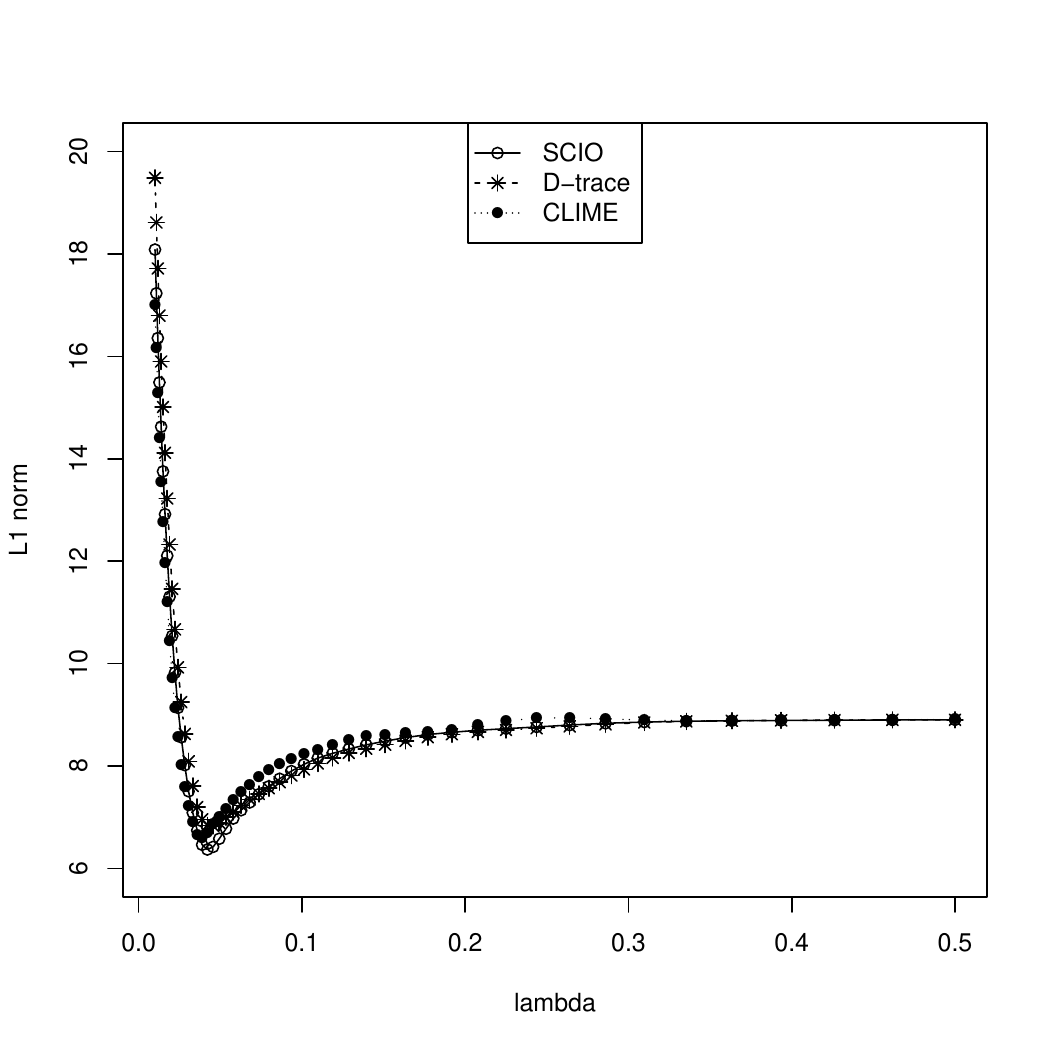, width=0.3\linewidth,angle=0}  \\
			(g) Frobenius norm for $\rho=0.8$  & (h) Spectral norm for  $\rho=0.8$ & (i) $L_1$ norm for  $\rho=0.8$\\
		\end{tabular}
	}
	\caption{Plots of the estimation errors versus the penalty parameter $\lambda$ under three sparsity levels based on the sample covariance matrix.}
	\label{fig:2}
\end{figure}

Figure \ref{fig:2} reports the performance of SCIO for three different sparsity levels $\rho=0.2$, $\rho=0.5$ and $\rho=0.8$. We plot the errors for the solution path with a series of tuning parameters and three methods: SCIO, D-trace and CLIME. From Figure \ref{fig:2}, we can see that these methods present similar patterns under all three norms. In other words, this demonstrates that SCIO performs similar as CLIME  which has been proved to be consistent under the $\ell_q$ sparsity condition.

\subsection{Robust matrix estimation}
The sub-Gaussian assumption is crucial in the analysis of the sample covariance matrix. To relax the assumption of exponential-type tails on the covariates,  \cite{avella2018robust} introduced a robust matrix estimator which only required a bounded fourth moment assumption. They constructed a Huber-type estimator for the population covariance matrix and got a robust estimator for the precision matrix by plugging the Huber-type estimator into the adaptively CLIME procedure  \citep{cai2016estimating}. 

Given the i.i.d. samples $\mX_k \in \mR^p$, $k=1,\ldots,n$,  \cite{avella2018robust} proposed to estimate the  covariance and the population mean based on the Huber loss function. In details, Huber's mean estimator $\tilde{\mu}_i$  satisfies the equation
\begin{align*}
	\sum_{k=1}^{n} \psi_{H}\left(\mX_{ki}-\tilde{\mu}_i\right)=0,
\end{align*} 
and the covariance estimator $\tilde{\sigma}_{ij}$ is defined by the equation
\begin{align*}
	\sum_{k=1}^{n} \psi_{H}\left(\mX_{ki}\mX_{kj}-(\tilde{\sigma}_{ij}+\tilde{\mu}_i\tilde{\mu}_j )\right)=0,
\end{align*} 
where $\psi_{H}(x)=\min \{H, \max (-H, x)\}$ denotes the Huber function.  Accordingly, we construct a robust estimator $\tilde{\boldsymbol{\Sigma}}=(\tilde{\sigma}_{ij})_{p \times p}$ and further project $\tilde{\boldsymbol{\Sigma}}$ to a cone of positive definite matrix
\begin{align*}
	\hbSig_2=\underset{\boldsymbol{\Sigma}\succeq\varepsilon\mathbf{I}}{\operatorname{argmin}}\|\boldsymbol{\Sigma}-\tilde{\boldsymbol{\Sigma}}\|_{\infty},
\end{align*}
where $\varepsilon$ is a small positive number.  This projection step can be easily implemented by the ADMM algorithm and see \cite{datta2017} for more details. 

\cite{avella2018robust} proposed to use $\hbSig_2$ as a pilot estimator and implemented the adaptively CLIME procedure to estimate the precision matrix. Here we study the SCIO method based on the robust estimator $\hbSig_2$. Following \cite{avella2018robust},  a bounded fourth moment condition is needed:
\bigskip

\noindent(A4).
Suppose that $\log{p}=o(n)$, and there exists a positive number $K>0$ such that
\begin{align*}
	\mE(X_i-\mathbf{\mu}_i)^4\leq K
\end{align*}
for all $1\leq i\leq p$.
\bigskip

Compared to the polynomial-type tails assumption (A3), the assumption (A4) here refines the moment order requirement from $4\gamma+4+\delta$ to $4$ and allows the result holds for a potentially larger $p$. The following proposition is from \cite{avella2018robust}.
\begin{proposition}[Proposition 3, \citealt{avella2018robust}]
	Under the assumption (A4),  for a sufficiently large constant $C$, we have
	\begin{align*}
		\mP\left(\|\hbSig_2-\bSig\|_\infty\geq  C\sqrt{\frac{\log{p}}{n}}\right)= O(p^{-\tau})
	\end{align*}
	for some constant $\tau>0$.
	\label{prop:2}
\end{proposition}

Similar to \cite{avella2018robust}, we can get an estimator $\tbOme_2$ by plugging  $\hbSig_2$ into SCIO and derive the convergence rates under different matrix norms.
\begin{corollary}
	Let $\lambda=C_0M_p\sqrt{\frac{\log p}{n}}$, where $C_0$ is a sufficiently large constant and $H=K(n / \log p)^{1 / 2}$ where $K$ is a given constant. For $\bOme\in\mathcal{U}_q(s_p,M_p)$, under assumptions (A1) and (A4), there exist sufficiently large constants $C_1,C_2$ satisfying that 
	\begin{align*}
		&\|\tbOme_2-\bOme\|_{\infty} \leq C_1 M_{p}^2 \sqrt{\frac{\log p}{n}},\\
		&\|\tbOme_2-\bOme\|_{L_1} \leq C_2 s_p M_{p}^{2-2q}\left(\frac{\log p}{n}\right)^{\frac{1}{2}(1-q)},
	\end{align*}
	with probability greater than $1-O(p^{-\tau})$, $\tau>0$.
	\label{cor:2}
\end{corollary}

\begin{remark}
	Note that \cite{avella2018robust} provided the optimal convergence rate with an additional technique assumption that the truncated population covariance matrix $\bSig_{H}=E\left\{1_{(\left|X_{u} X_{v}\right|\leq H)} X_{u} X_{v}\right\}$ satisfies that $\|\bSig_{H}\bOme-\mathbf{I}\|_{\infty}\leq C\sqrt{\frac{\log{p}}{n}}$. Although the convergence rate provided in Corollary \ref{cor:2} is not optimal, the optimal rate can be readily obtained by imposing the same condition on the truncated population covariance matrix in \cite{avella2018robust}.
\end{remark}

\begin{figure}
	\centerline{
		\begin{tabular}{ccc}
			\psfig{figure=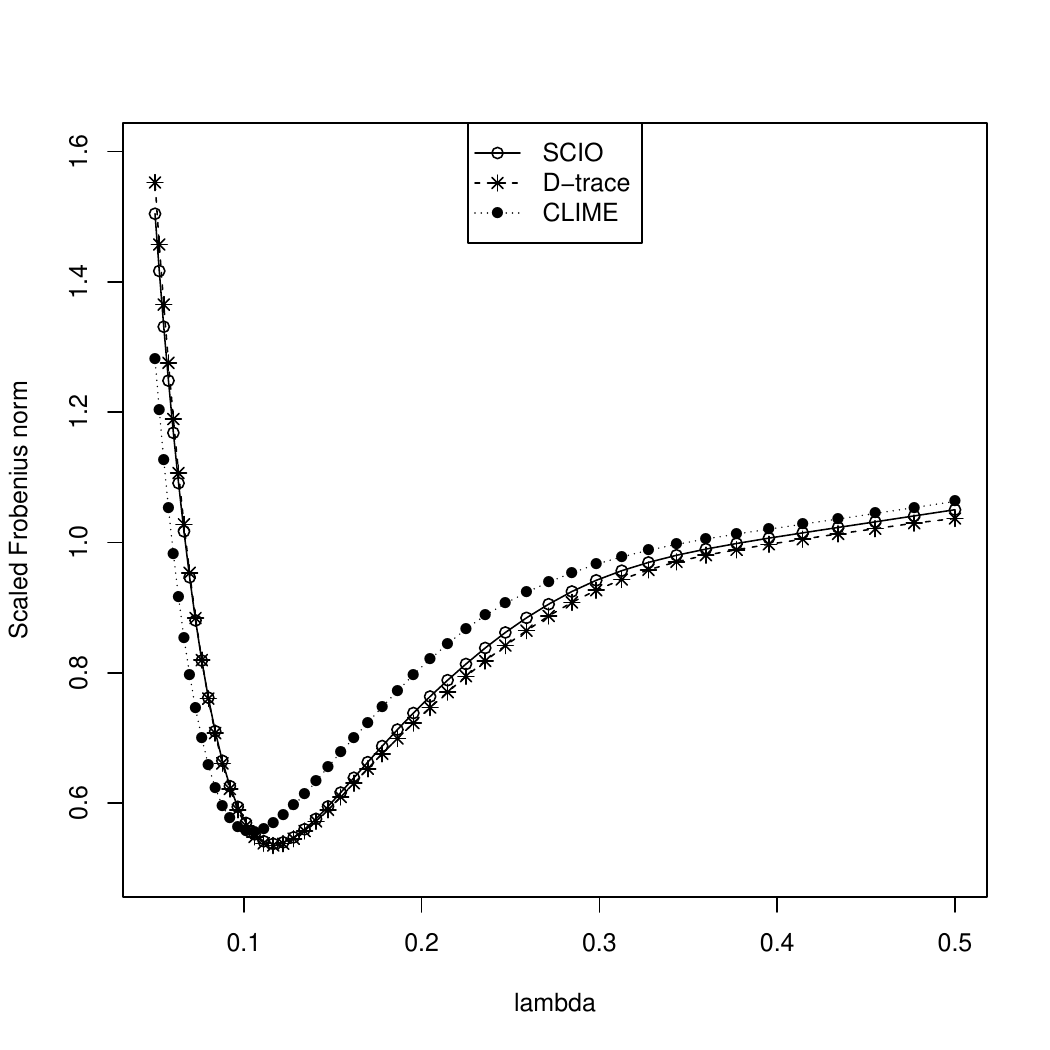, width=0.3\linewidth,angle=0} &
			\psfig{figure=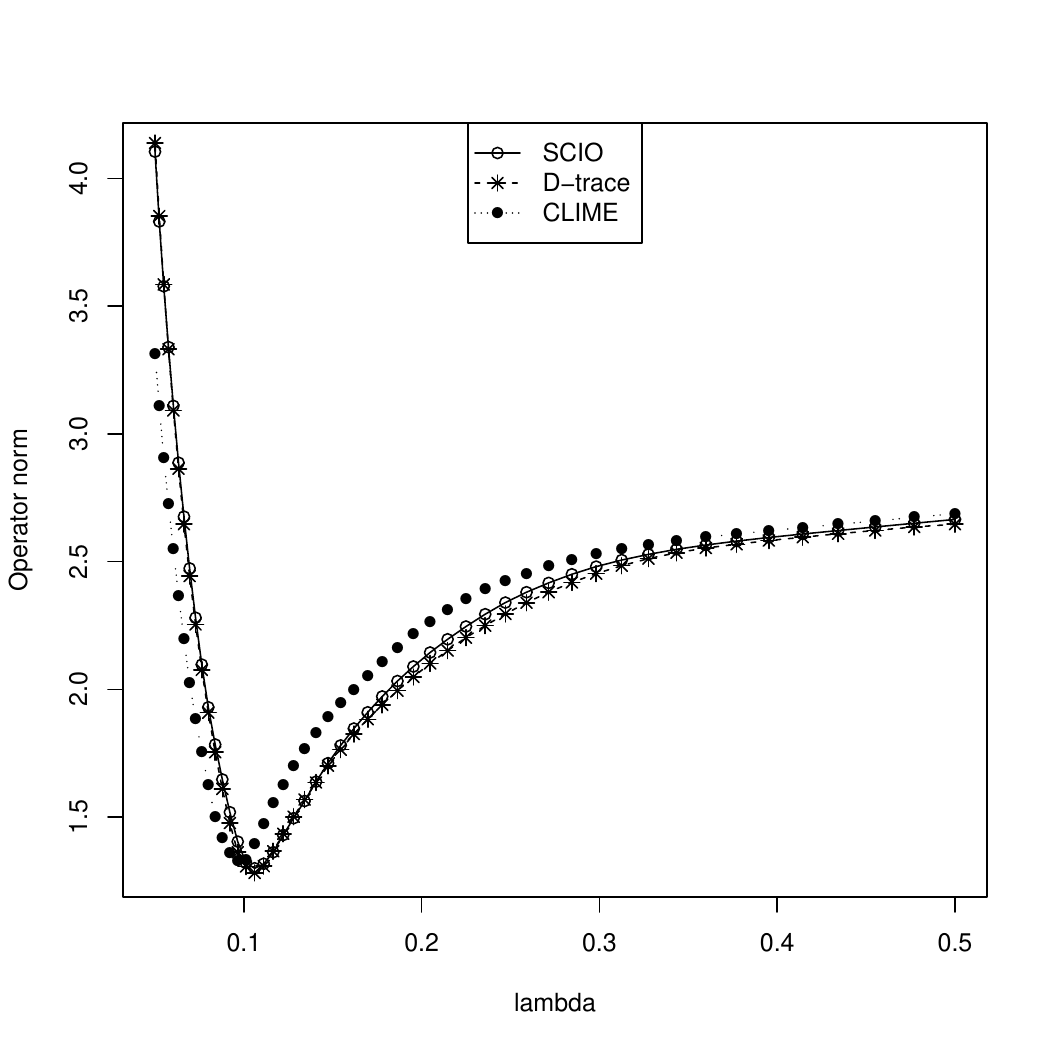, width=0.3\linewidth,angle=0}&
			\psfig{figure=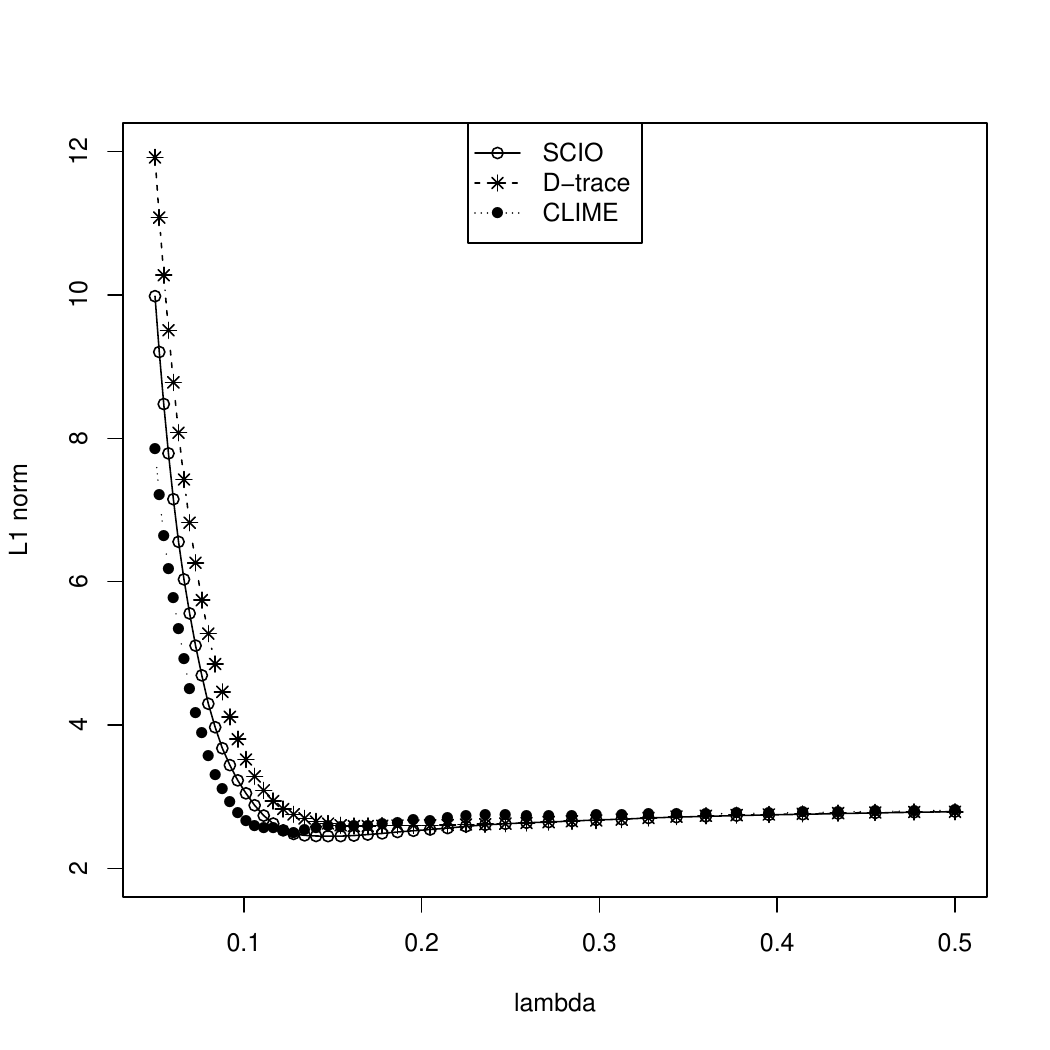, width=0.3\linewidth,angle=0}  \\
			(a) Frobenius norm & (b) Spectral norm & (c) $L_1$ norm
		\end{tabular}
	}
	\caption{Plots of the estimation errors versus the penalty parameter $\lambda$ based on the Huber-type estimator. }
	\label{fig:3}
\end{figure}

To demonstrate these results numerically, we repeat the second scenario in \cite{avella2018robust}  where the data $\mathbf{X}$ is generated from a student $t$ distribution with 3.5 degrees of freedom and infinite kurtosis.  Here the sub-Gaussian assumption is void.  We still consider the precision matrix $\bOme=(0.5^{|i-j|})_{p \times p}$. The sample size $n$ is set to 200 and the dimension $p$ is 100. Figure \ref{fig:3} reports the numeric performances of SCIO, D-trace and CLIME based on the robust estimator $\hbSig_2$.  All three methods perform comparably and align well for the heavy-tailed distribution.

Furthermore, we conduct a numerical simulation to illustrate the robustness of our SCIO estimator for the heavy-tailed distribution. In details, we plug the sample covariance matrix $\hbSig_1$ and the Huber-type covariance matrix estimator $\hbSig_2$ into the SCIO.  We set the sample size as  $n=100$ and generate the data matrix $\mathbf{X}$ with a multivariate $t$ distribution with 5 degrees of freedom, zero mean and a covariance matrix $\bSig=\bOme^{-1}$, where $\bOme=(\rho^{|i-j|})_{p \times p}$ and $\rho=0.2,0.5$. Note that this distribution only has $4$-th order moment.  Table \ref{tab3} reports the spectral norm error for different dimensions $p$ based on 50 replications. From Table \ref{tab3},  we can see that the Huber-type precision matrix estimator performs better than the one with the sample covariance matrix.  This result is consistent with our theoretical improvement from the requirement of Assumption (A3) to a milder one (A4).

\begin{table*}[htbp]

	\centering
	
	\caption{Comparison of SCIO with the sample and the Huber-type covariance matrices for the heavy-tailed data.}
	
	\label{tab3}
	
	\begin{tabular}{cccccccccc}
		
		\toprule
		\toprule
		
		\multirow{2}{*}{$\emph{p}$} & \multicolumn{2}{c}{$\rho=0.2$} & \multicolumn{2}{c}{$\rho=0.5$} \\
		
		\cmidrule(r){2-3} \cmidrule(r){4-5} 
		
		&  $\hbSig_1$     &  $\hbSig_2$
		
		&  $\hbSig_1$     &  $\hbSig_2$\\
		
		\midrule
		$100 $             &0.75(0.02)                         & 0.67(0.05)                    & 1.69(0.02)                   & 1.54(0.01)           \\
		
		$200$             &0.78(0.01)                          & 0.68(0.01)                    & 1.73(0.10)                   & 1.63(0.02)           \\
		
		$400$             &0.85(0.02)                          & 0.70(0.07)                   & 2.00(0.15)                   & 1.70(0.04)          \\
		
		$600$             &0.88(0.05)                          & 0.74(0.05)                   & 2.17(0.11)                   & 1.85(0.07)           \\
		
		\bottomrule
		
	\end{tabular}
	
\end{table*}

\subsection{Non-parametric rank-based estimation}
For Gaussian distributions, the precision matrix characterizes the conditional independence among covariates. For non-Gaussian data, \cite{liu2009nonparanormal} introduced a non-paranormal graphical model.  \cite{liu2012high} and
\cite{xue2012regularized} studied the precision matrix estimation for this non-paranormal graphical model where 
the precision matrix: $\bOme={\bSig}^{-1}$ was defined by the transformed samples and $\bSig$ was the correlation matrix. In details, they proposed to estimate the correlation matrix by the non-parametric rank-based statistics such as Spearman's rho and Kendall's tau.

Given the sample data matrix $(X_{ij})_{n \times p}=(\mX_1,\cdots,\mX_n)\trans$, we convert them to rank statistics denoted by $(r_{ij})_{n \times p}=(\r_1,\cdots,\r_n)\trans$ where each column $(r_{1j},\cdots,r_{nj})$ serves as the rank statistic of $(X_{1j},\cdots,X_{nj})$.  Spearman's rho correlation coefficient $\widehat{\rho}_{i j}$ is defined as the Pearson correlation between the columns $\r_i$ and $\r_j$, that is, 
\begin{align*}
	\mbox{Spearman's rho:}~ \widehat{\rho}_{i j}=\frac{\sum_{k=1}^{n}\left(r_{ki}-\bar{r}_{i}\right)\left(r_{kj}-\bar{r}_{j}\right)}{\sqrt{\sum_{k=1}^{n}\left(r_{ki}-\bar{r}_{i}\right)^2 \cdot\sum_{k=1}^{n}\left(r_{kj}-\bar{r}_{j}\right)^2}}. 
\end{align*}	
Similarly, Kendall's tau correlation coefficient is defined by
\begin{align*}	
	\mbox{Kendall's tau:}~  \widehat{\tau}_{ij}=\frac{2}{n(n-1)} \sum_{1 \leq k<k^{\prime} \leq n} \operatorname{sign}\{(X_{ki}-X_{k^{\prime}i})(X_{kj}-X_{k^{\prime}j})\}.
\end{align*}
Based on Spearman's rho and Kendall's tau correlation coefficients, we are able to construct two non-parametric  estimators $\tilde{\bSig}_{3\rho}$ and $\tilde{\bSig}_{3\tau}$ for the correlation matrix $\bSig$, where 
\begin{align*}
	(\tilde{\bSig}_{3\rho})_{ij}=\left\{\begin{array}{ll}
		2 \sin \left(\frac{\pi}{6} \widehat{\rho}_{ij}\right), & i \neq j \\
		1, & i=j
	\end{array}\right.
\end{align*}
and
\begin{align*}
	(\tilde{\bSig}_{3\tau})_{ij}=\left\{\begin{array}{ll}
		\sin \left(\frac{\pi}{2} \widehat{\tau}_{ij}\right), & i \neq j \\
		1, & i=j
	\end{array}\right..
\end{align*}
Moreover, we still need an additional projection step
\begin{align*}
	\hbSig_3=\argmin_{\bSig \succeq \varepsilon\mathbf{I}} \|\bSig-\tilde{\bSig}_3\|_{\infty}
\end{align*}
to obtain the final positive definite estimator $\hbSig_{3\rho}$ or $\hbSig_{3\tau}$. Note that if $\mathbf{X}$ satisfies the non-paranormal distribution,  \cite{liu2012high} proved that $\tilde{\bSig}_{3\rho}$ and $\tilde{\bSig}_{3\tau}$ are consistent estimators of $\bSig$ under the element-wise $\ell_\infty$ norm.   The following proposition is from \cite{liu2012high}.
\begin{proposition}[Theorem 4.1 and 4.2, \citealt{liu2012high}]
	Assuming that $\mathbf{X}$ satisfies a non-paranormal distribution, there exists a sufficiently large constant $C$ such that
	\begin{align*}
		&\mP\left(\|\hbSig_{3\rho}-\bSig\|_\infty\geq  C\sqrt{\frac{\log{p}}{n}}\right)= O(p^{-1}),\\
		&\mP\left(\|\hbSig_{3\tau}-\bSig\|_\infty\geq  C\sqrt{\frac{\log{p}}{n}}\right)= O(p^{-1}).
	\end{align*}
	\label{prop:3}
\end{proposition}

To estimate the sparse precision matrix, 
\cite{liu2012high} proposed to plug $\hbSig_{3\rho}$ or $\hbSig_{3\tau}$ into the graphical Dantzig selector \citep{yuan2010high}, CLIME\citep{Cai2011A}, the graphical Lasso  \citep{friedman2008sparse},  or the neighborhood pursuit estimator \citep{meinshausen2006high}. In this part, we consider the SCIO procedure with $\hbSig_{3\rho}$ or $\hbSig_{3\tau}$. Denote $\tbOme_3$ as the precision matrix estimator by plugging $\hbSig_{3\rho}$ or $\hbSig_{3\tau}$ into SCIO. The following corollary holds for both $\hbSig_{3\rho}$ and $\hbSig_{3\tau}$.
\begin{corollary}
	Let $\lambda=C_0M_p\sqrt{\frac{\log p}{n}}$, where $C_0$ is a sufficiently large constant. For $\bOme\in\mathcal{U}_q(s_p,M_p)$, under assumptions (A1) and that $\mathbf{X}$ satisfies a non-paranormal distribution, there exist sufficiently large constants $C_1,C_2$ satisfying that 
	\begin{align*}
		&\|\tbOme_3-\bOme\|_{\infty} \leq C_1 M_{p}^2 \sqrt{\frac{\log p}{n}},\\
		&\|\tbOme_3-\bOme\|_{L_1} \leq C_2 s_p M_{p}^{2-2q}\left(\frac{\log p}{n}\right)^{\frac{1}{2}(1-q)}
	\end{align*}
	with probability greater than $1-O(p^{-1})$.
	\label{cor:3}
\end{corollary}

\begin{figure}[]
	\centerline{
		\begin{tabular}{ccc}				
			\psfig{figure=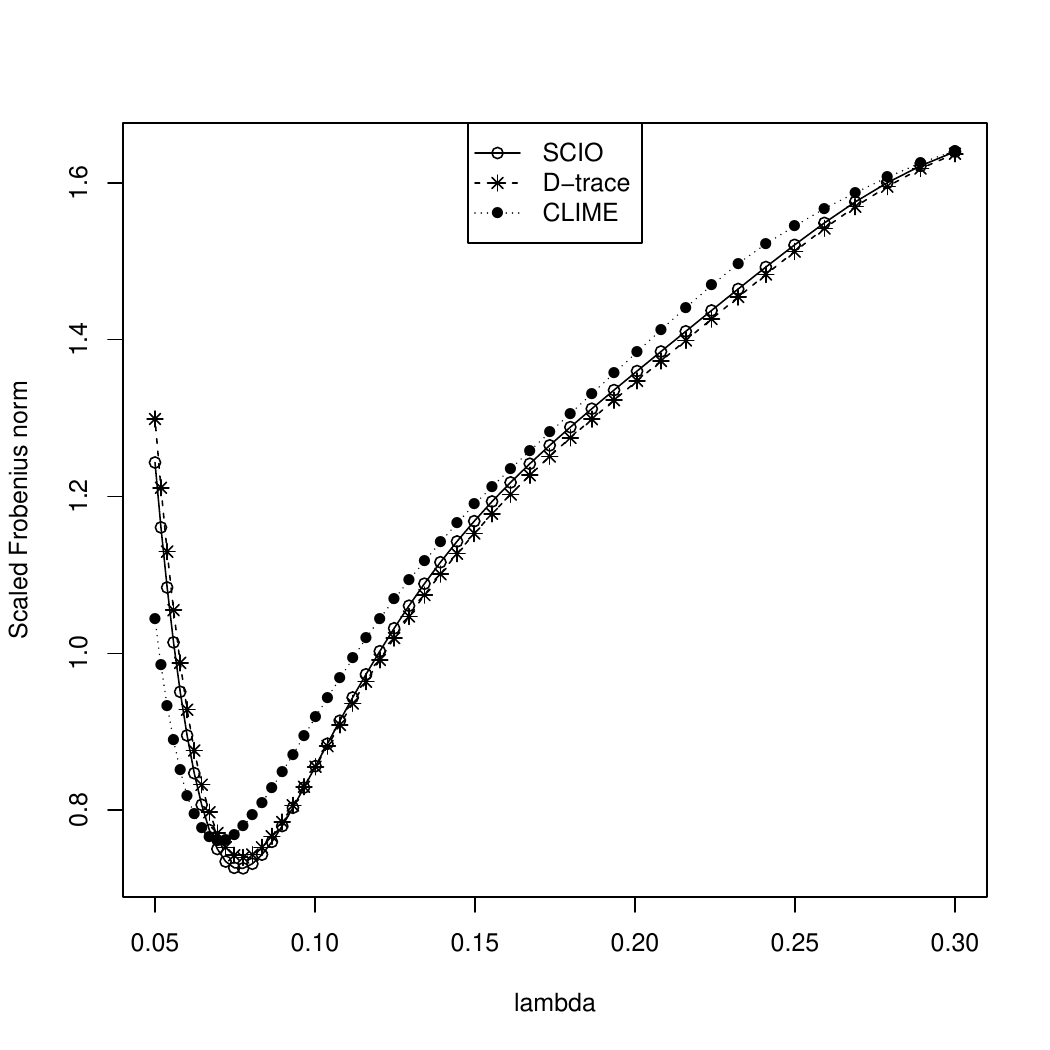, width=0.3\linewidth,angle=0} &
			\psfig{figure=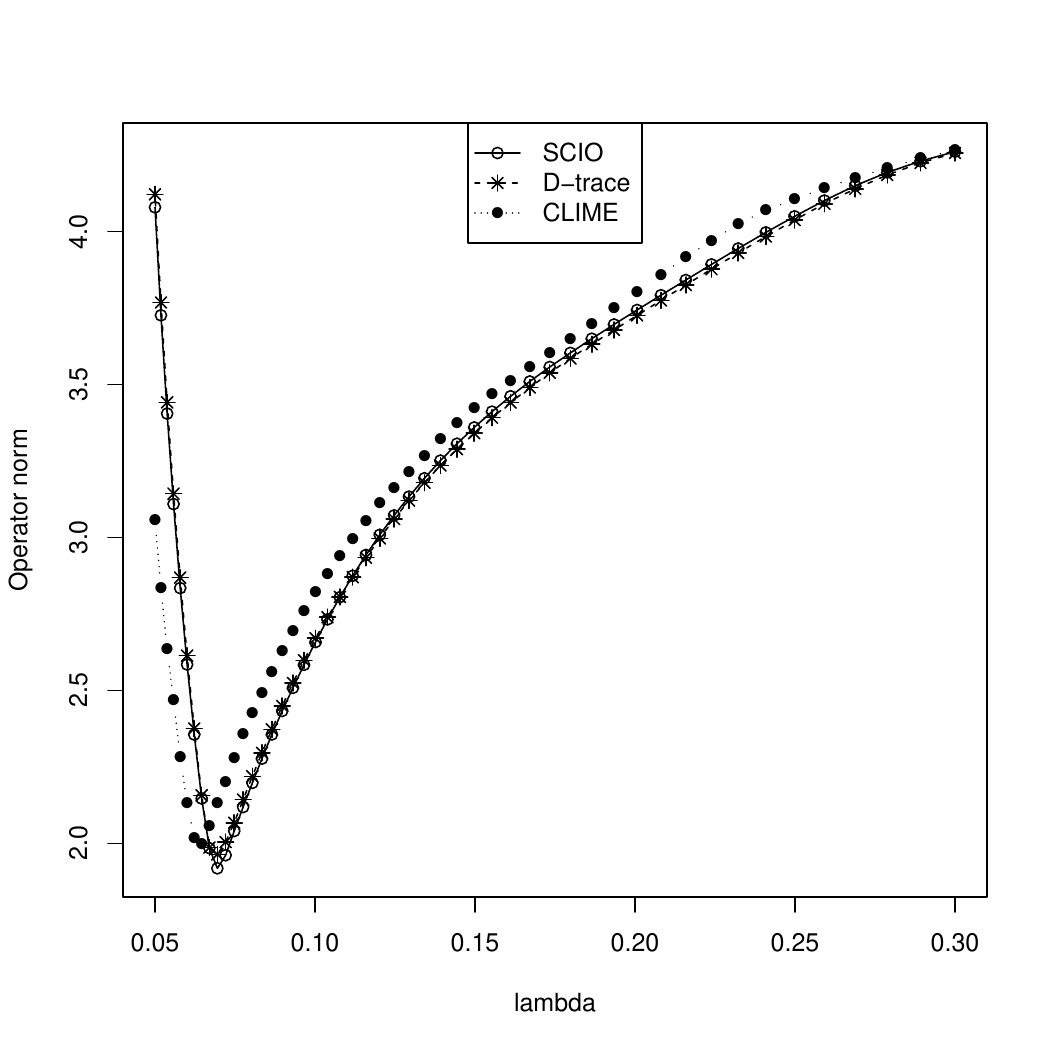, width=0.3\linewidth,angle=0}&
			\psfig{figure=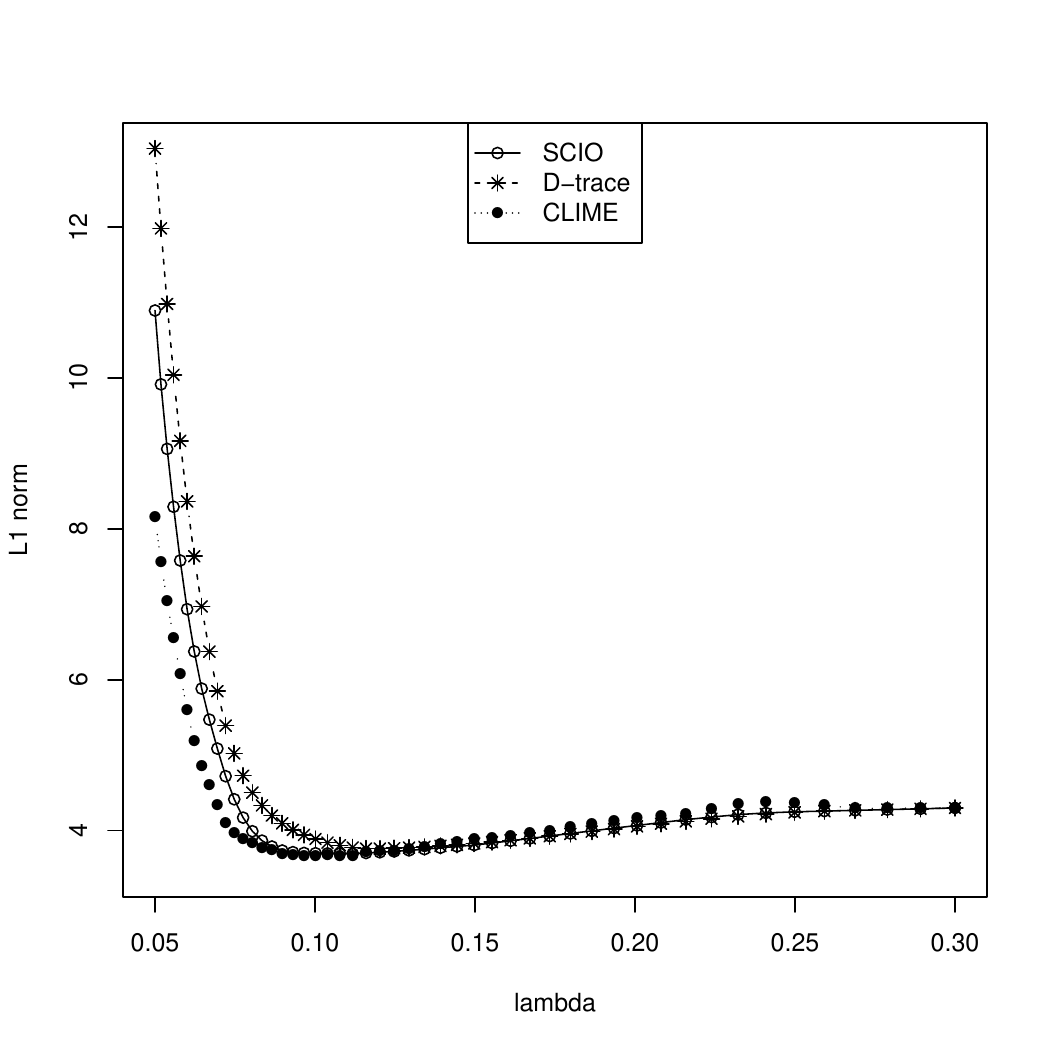, width=0.3\linewidth,angle=0}  \\
			(a) Frobenius norm & (b) Spectral norm & (c) $L_1$ norm \\
			\psfig{figure=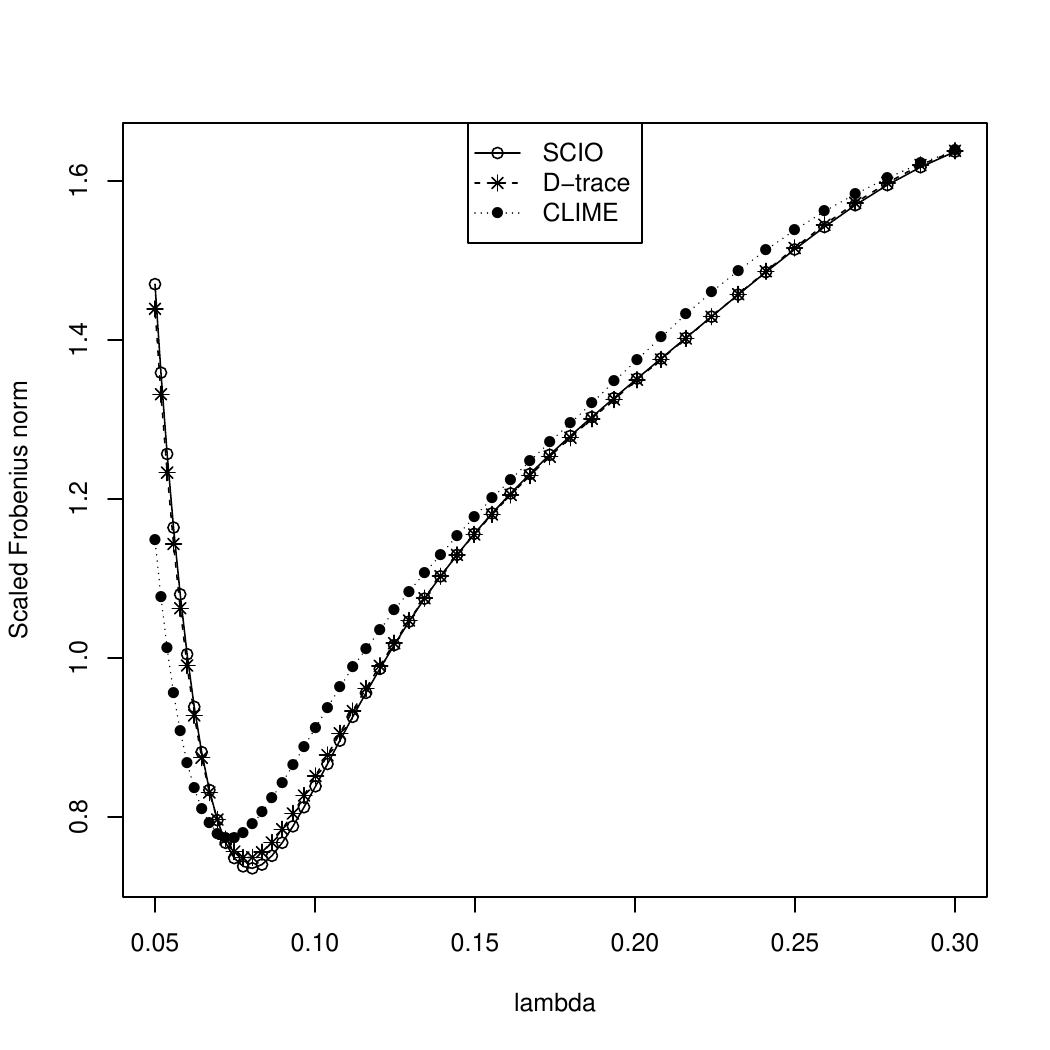, width=0.3\linewidth,angle=0} &
			\psfig{figure=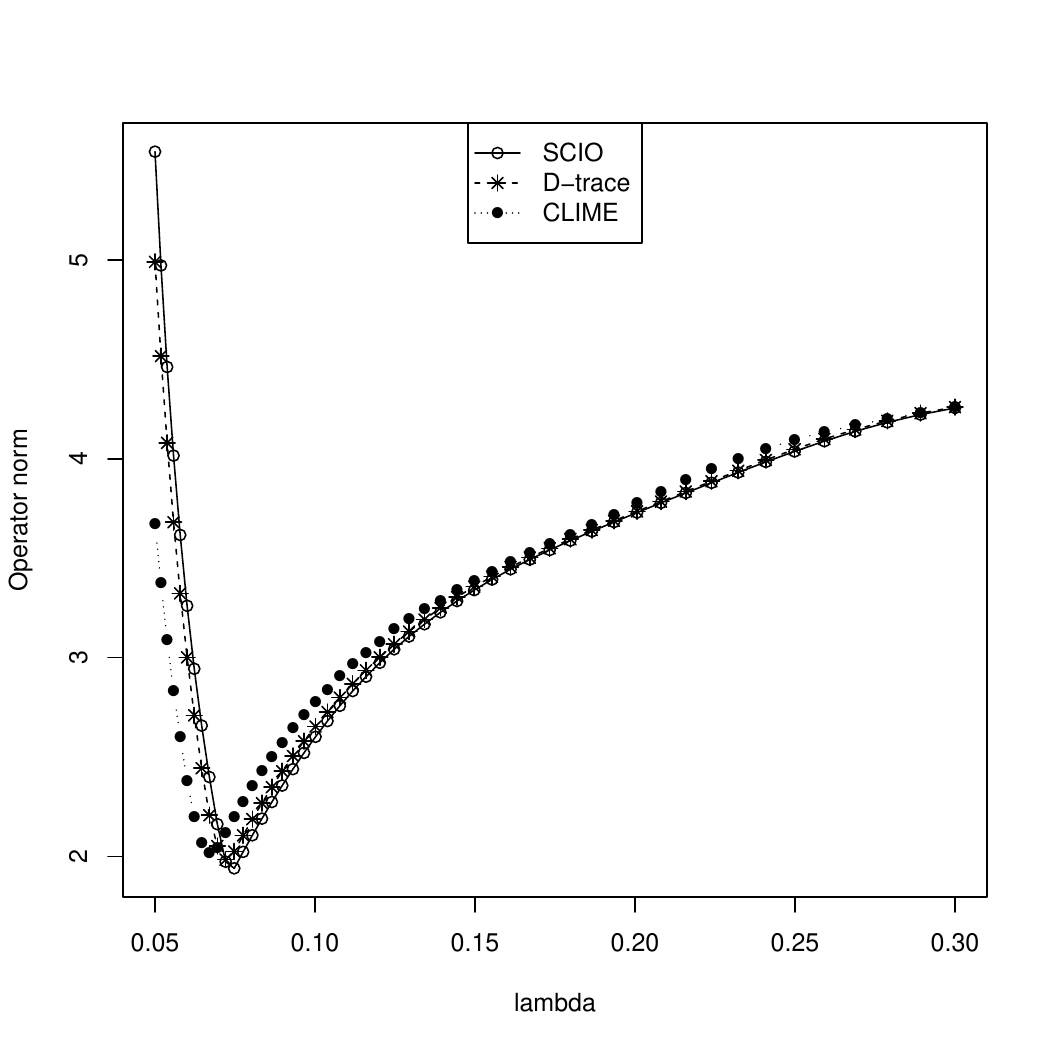, width=0.3\linewidth,angle=0}&
			\psfig{figure=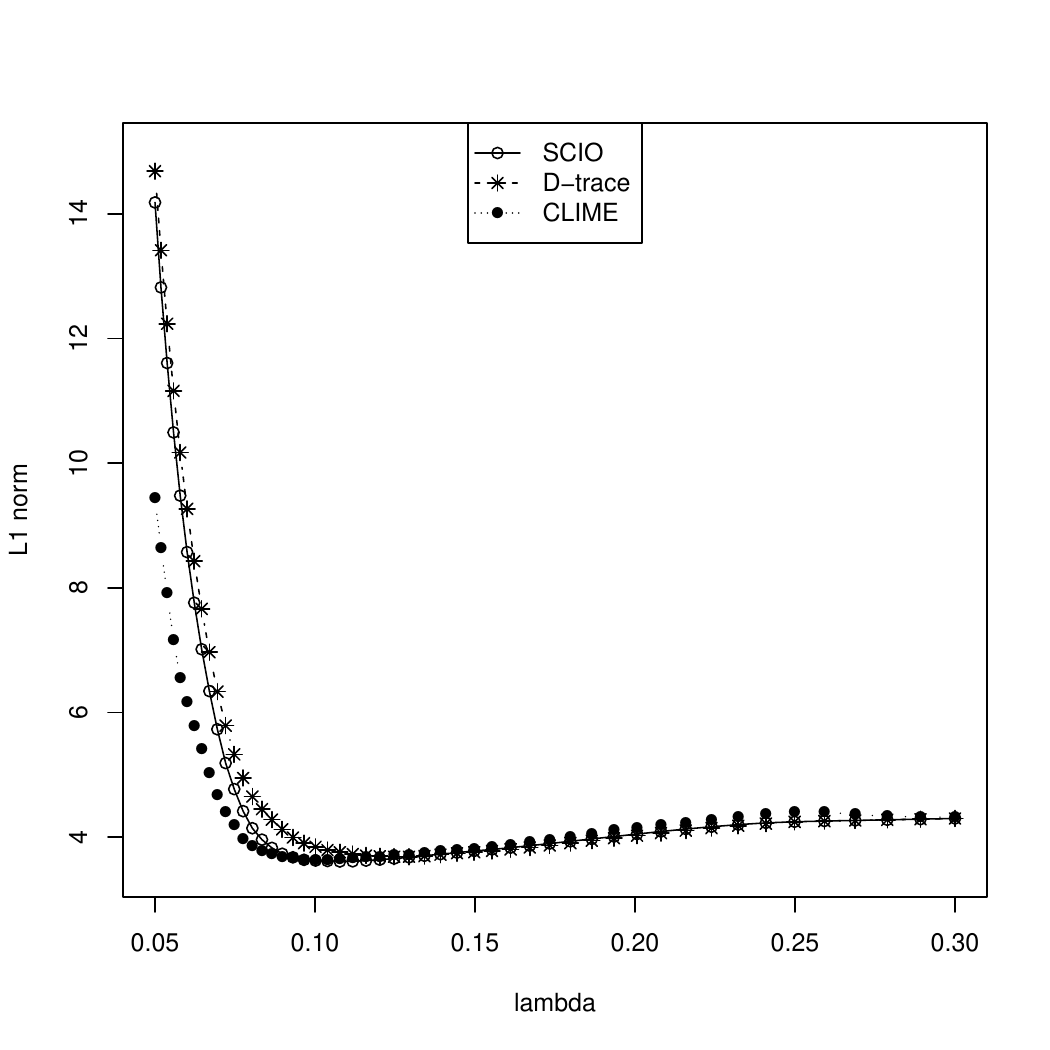, width=0.3\linewidth,angle=0}  \\
			(d) Frobenius norm& (e) Spectral norm& (f) $L_1$ norm\\	
		\end{tabular}
	}
	\caption{Plots of the estimation errors versus the penalty parameter $\lambda$ based on Spearman's rho and Kendall's tau estimation. Here (a)-(c) are the results for Spearman estimation and (d)-(f) are the results for Kendall estimation.}
	\label{fig:4}
\end{figure}

To conduct numeric simulations, we assume that $\mathbf{X}$ follows a non-paranormal distribution $f(\mathbf{X})\sim\mathcal{N}(0,\boldsymbol{\Sigma})$. Following Definition 9 in \cite{liu2009nonparanormal} and Definition 5.1 in \cite{liu2012high}, we choose the transformation function $f$ as the Gaussian CDF transformation function with $\mu_{g_0}=0.05$ and $\sigma_{g_0}=0.4$. To mimic the weak sparse case, we consider $\bSig$ as the correlation matrix of $\bOme_0^{-1}$, where $\bOme_0=(0.5^{|i-j|})_{p \times p}$. We set the sample size $n=200$ and the dimension $p=100$ again. Based on Spearman's rho estimation $\hbSig_{3\rho}$ or Kendall's tau estimation $\hbSig_{3\tau}$, Figure \ref{fig:4} reports the numeric performances of the SCIO, D-trace and CLIME.  Again,  we can see that all three methods perform comparably.

As \cite{avella2018robust} showed, Proposition \ref{prop:3} works for the elliptically distributed $\mathbf{X}$, which includes the multivariate $t$ distributed random variables. Here, we evaluate the robust performance of non-parametric rank-based SCIO estimation under the heavy-tailed circumstance. The setting is the same as the one of Table  \ref{tab3}. Table \ref{tab4} shows the spectral norm error of SCIO with Pearson's correlation matrix $\hbSig_1$, Spearman's rho correlation matrix $\hbSig_{3\rho}$ and Kendall's tau correlation matrix $\hbSig_{3\tau}$. We can see that SCIO with  non-parametric correlation estimators outperform the one with the Pearson's correlation matrix. Under our settings, Spearman's estimation performs slightly better than Kendall's estimation. The numerical results verify the robustness of our non-parametric rank-based SCIO estimation for the heavy-tailed case.

\begin{table*}[htbp]
	
	\scriptsize
	
	\centering
	
	\caption{Comparison of precision matrix estimation errors under the spectral norm for the non-parametric estimators over 50 replications.}
	
	\label{tab4}
	
	\begin{tabular}{cccccccccc}
		
		\toprule
		\toprule
		
		\multirow{2}{*}{$\emph{p}$} & \multicolumn{3}{c}{$\rho=0.2$} & \multicolumn{3}{c}{$\rho=0.5$} \\
		
		\cmidrule(r){2-4} \cmidrule(r){5-7} 
		
		&  $\hbSig_1$     &  $\hbSig_{3\rho}$
		
		&  $\hbSig_{3\tau}$      &  $\hbSig_1$
		
		&  $\hbSig_{3\rho}$      &  $\hbSig_{3\tau}$     \\
		
		\midrule
		$100 $             &0.79(0.08)                          & 0.72(0.01)                    & 0.73(0.01)                   & 3.01(0.58)           & 2.78(0.06)           & 2.99(0.14)        \\
		
		$200$             &0.79(0.04)                          & 0.73(0.01)                    & 0.75(0.01)                   & 3.31(0.41)           & 3.08(0.07)           & 3.20(0.11)        \\
		
		$400$             &0.81(0.03)                          & 0.75(0.01)                   & 0.77(0.02)                   & 3.51(0.13)          & 3.45(0.11)           & 3.47(0.13)         \\
		
		$600$             &0.86(0.04)                          & 0.77(0.02)                   & 0.78(0.02)                   & 3.63(0.10)           & 3.53(0.12)          & 3.54(0.12)         \\
		
		\bottomrule
		
	\end{tabular}
	
\end{table*}

\subsection{Matrix data estimation}
The matrix variate data are frequently encountered in real applications where the covariance matrix has a Kronecker product structure $\bSig=\A\otimes\B$. To study the matrix data, it is of great interest to estimate the graphical structures $\bOme=\bSig^{-1}=(\A\otimes\B)^{-1}$ \citep{leng2012sparse, zhou2014gemini}. For the brevity, we assume $\A$ and $\B$ are all correlation matrices, which means the diagonal entries are ones.  

Given the i.i.d. matrix samples $\mX(t) \in \mR^{f \times m},~t=1,\ldots,n$,  \cite{zhou2014gemini} developed the Gemini estimator for the precision matrix $\A^{-1}\otimes \B^{-1}$. Writing the columns of $\mX(t)$ as $x(t)^{1},\cdots,x(t)^m \in \mR^f$,  \cite{zhou2014gemini} proposed to estimate $\A$ by  
\begin{align*}
	(\hbSig_{4A})_{ij}=\frac{\sum_{t=1}^n (x(t)^{i})\trans (x(t)^{j})  }{\sqrt{\sum_{t=1}^n (x(t)^{i})\trans (x(t)^{i}) } \sqrt{\sum_{t=1}^n (x(t)^{j})\trans (x(t)^{j})}}
\end{align*}
and the estimator $\hbSig_{4B}$ for $\B$ is constructed similarly based on the rows of $\mX(t)$.  The final estimation of $\bOme$ is obtained by implementing the graphical Lasso  or CLIME with $\hbSig_{4A}$ and $\hbSig_{4B}$. \cite{zhou2014gemini} derived the convergence rate under $\ell_0$ sparsity condition for the graphical Lasso  and introduced the CLIME procedure to refine their convergence rates. For the $\ell_q$ sparse matrix, \cite{zhou2014gemini} did not provide the explicit theoretical results.

In this part, we study the SCIO method based on $\hbSig_{4A}$ and $\hbSig_{4B}$.  We first present an approximate sparsity condition for $\A^{-1}$ and $\B^{-1}$.

\bigskip
\noindent(A5). Suppose $\A^{-1} \in\mathcal{U}_q(s_m,M_m)$ and $\B^{-1} \in\mathcal{U}_q( \tilde{s}_f, \tilde{M}_f)$ for a given $q \in [0,1)$. Moreover, the  parameters $s_m,  \tilde{s}_f, M_m,  \tilde{M}_f$ satisfy 
\begin{align*}
	&s_mM_m^{2-q}=o\left(\frac{\sqrt{nf}}{\log^{\frac{1}{2}}(m\vee f)}\right),\\ & \tilde{s}_f \tilde{M}_f^{2-q}=o\left(\frac{\sqrt{nm}}{\log^{\frac{1}{2}}(m\vee f)}\right).
\end{align*}
The following proposition is  from Theorem 4.1 of \cite{zhou2014gemini}.
\begin{proposition}[Theorem 4.1, \citealt{zhou2014gemini}] For $t=1,\ldots, n$, 
	suppose that $\operatorname{vec}(\mX(t)) \sim \mathcal{N}_{f, m}(0, \A \otimes \mathbf{B})$. Under the assumption (A5) and the assumption (A2) in \cite{zhou2014gemini}, there exists a sufficiently large constant $C$ such that
	\begin{align*}
		&\mP\left(\|\hbSig_{4A}-\A\|_\infty\geq  C\sqrt{\frac{\log(m\vee f)}{nf}} \right)= O((m \vee f)^{-2}),\\
		&\mP \left(\|\hbSig_{4B}-\B\|_\infty\geq  C\sqrt{\frac{\log(m\vee f)}{nm}}\right)= O((m \vee f)^{-2}).
	\end{align*}
	\label{prop:4}
\end{proposition}
As an application of our Theorem \ref{thm:2}, we can derive the theoretical result of the SCIO estimator for estimating $\bOme=\A^{-1}\otimes \mathbf{B}^{-1}$.
\begin{corollary}
	Suppose that $\operatorname{vec}(\mX(t)) \sim \mathcal{N}_{f, m}(0, \A_{m\times m} \otimes\mathbf{B}_{f\times f}),~t=1,\ldots, n$. Let $\lambda_A=C_AM_m\frac{\log^{\frac{1}{2}}(m\vee f)}{\sqrt{fn}}$ and $\lambda_B=C_B \tilde{M}_f\frac{\log^{\frac{1}{2}}(m\vee f)}{\sqrt{mn}}$, where $C_A$ and $C_B$ are sufficiently large constants. Under our assumption (A5) and the assumption (A2) in \cite{zhou2014gemini},  there exists a sufficiently large constant $C$ such that 
	\begin{align*}
		\|\tbOme_4-\bOme\|_{2} \leq C \left(\frac{\log(m\vee f)}{n}\right)^{\frac{1-q}{2}}  \left(s_mM_m^{2-2q}f^{(q-1)/2} + \tilde{s}_f \tilde{M}_f^{2-2q} m^{(q-1)/2}\right)
	\end{align*}
	with probability greater than $1-O((m \vee f)^{-2})$.
	\label{cor:4}
\end{corollary}
Compared with Theorem 3.3 of \cite{zhou2014gemini}, Corollary \ref{cor:4} is derived under  a general $\ell_q$ sparsity condition. In particular, the special case $q=0$ corresponds to Theorem 3.3 of \cite{zhou2014gemini} and these two results are consistent due to the dual properties between Lasso  and the Dantzig selector. Moreover, our result can be generalized to the sub-Gaussian condition of the matrix data \citep{hornstein2019joint} and we omit the details.

\begin{figure}
	\centerline{
		\begin{tabular}{ccc}				
			\psfig{figure=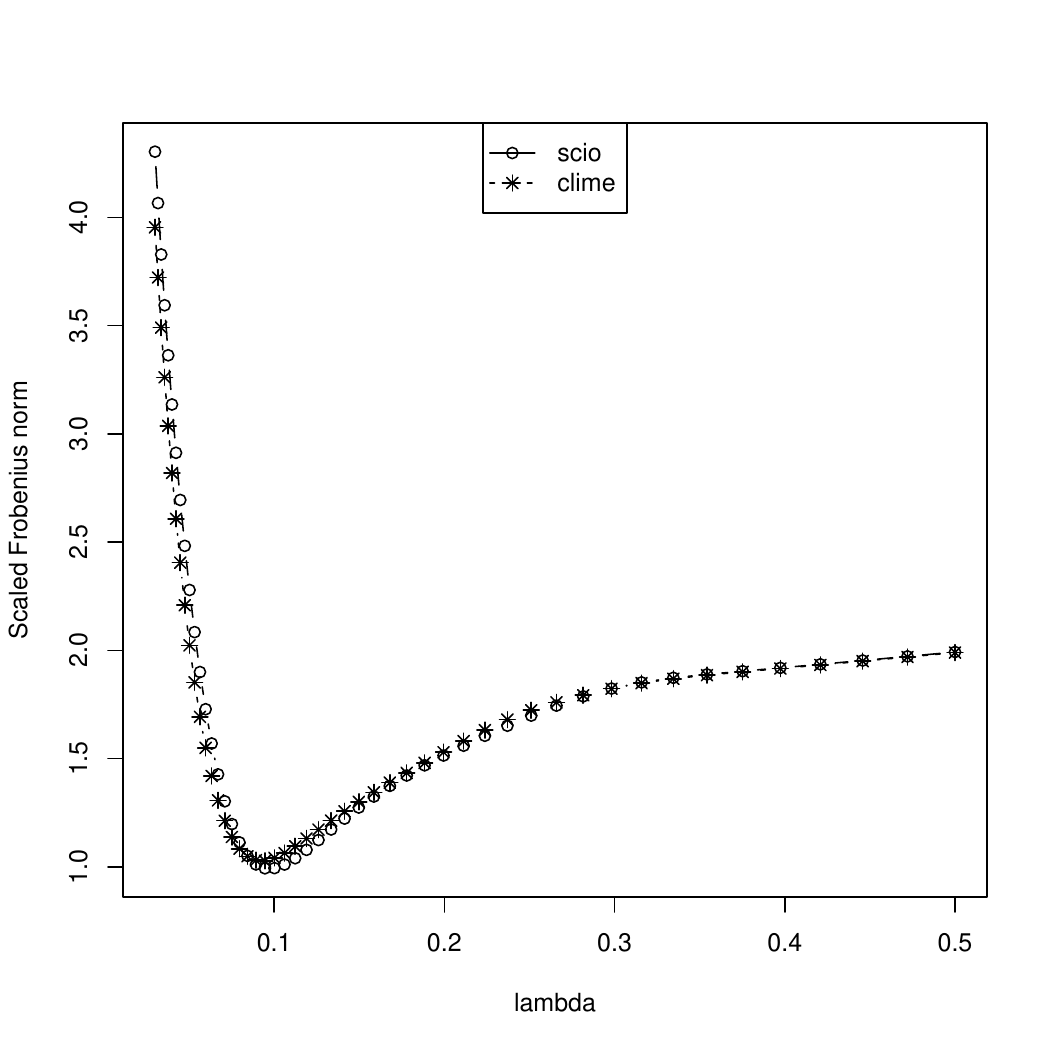, width=0.3\linewidth,angle=0} &
			\psfig{figure=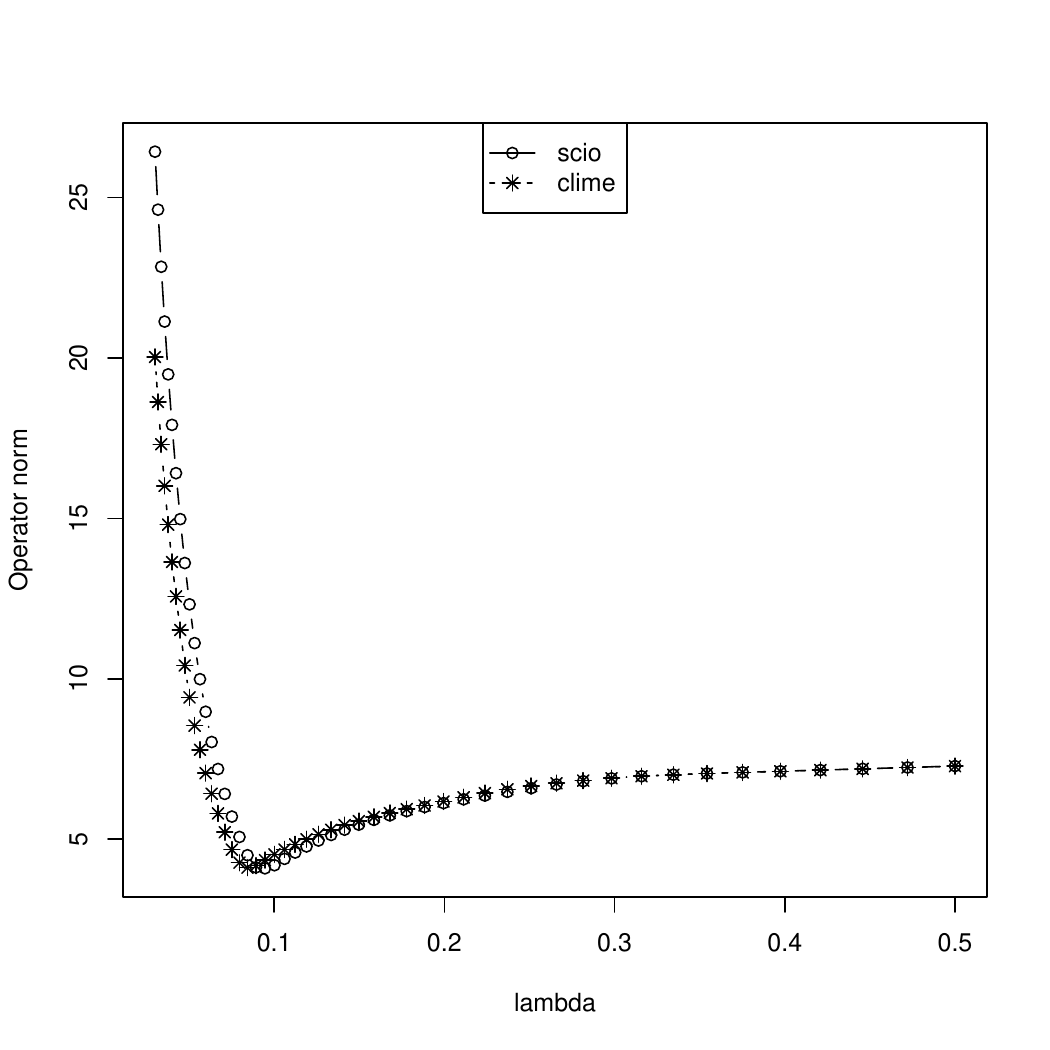, width=0.3\linewidth,angle=0}&
			\psfig{figure=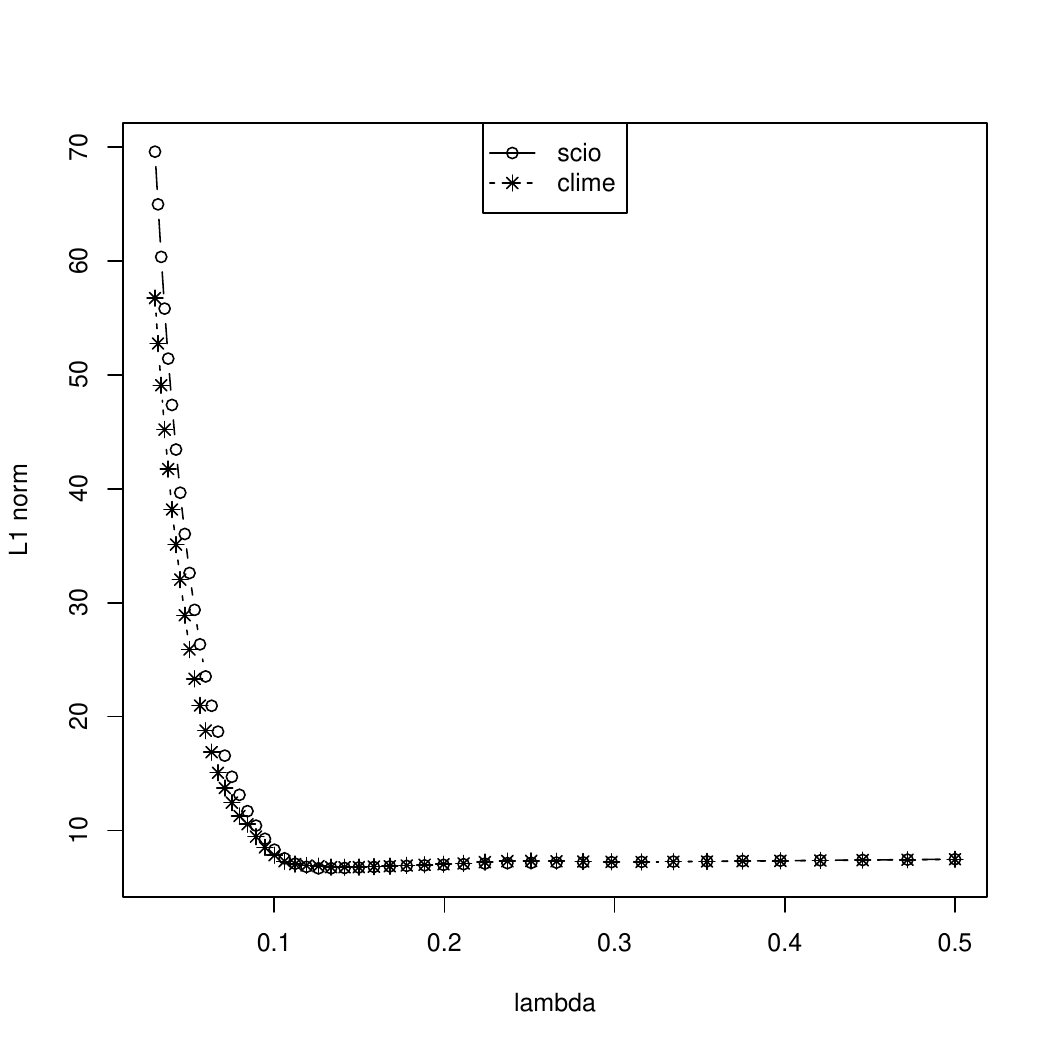, width=0.3\linewidth,angle=0}  \\
			(a) Frobenius norm & (b) Spectral norm & (c) $L_1$ norm \\
		\end{tabular}
	}
	\caption{Plots of the estimation errors versus the penalty parameter $\lambda$ based on the Gemini estimator  for the matrix data.}
	\label{fig:5}
\end{figure}

To conduct the simulations, we generate the data from the matrix normal distribution  $\operatorname{vec}(\mX(t)) \sim \mathcal{N}_{f, m}(0, \A_{m\times m}\otimes \mathbf{B}_{f\times f})$. To mimic the $\ell_q$ sparsity, we choose $\mathbf{B}_{ij}$ as the correlation matrix of $\boldsymbol{\Phi}^{-1}$, where $\boldsymbol{\Phi}_{ij}=(0.2^{|i-j|})_{f \times f}$ and $\mathbf{A}_{ij}$ as the correlation matrix of $\boldsymbol{\Theta}^{-1}$, where $\boldsymbol{\Theta}_{ij}=(0.5^{|i-j|})_{m \times m}$ . We set the dimension of $\mathbf{A}$ as 80 and the dimension of $\mathbf{B}$ as 40. The sample size $n$ is taken as 3. Figure \ref{fig:5} reports the performance of the Gemini estimator under several matrix norms where the penalty level $\lambda$ of $\mathbf{A}^{-1}$ is varying and the penalty level of $\mathbf{B}^{-1}$ is set to $0.15$ for simplicity. From Figure \ref{fig:5}, we can observe that the Gemini method based on SCIO performs similarly as the Gemini method based on CLIME, which means SCIO is also applicable to the matrix data.

\section{Discussion}

This article revisits the SCIO method proposed by \cite{liu2015fast} and explores the theoretical and numerical properties of SCIO under the weak sparsity condition. Intuitively, the approach to obtain our matrix estimation error bound by plugging in the sample covariance matrix is similar to the process of obtaining the prediction bound in regression setting. For the classical Lasso  problem, \cite{ye2010rate} and \cite{sun2012scaled} have analyzed the Lasso  method or its variants under general weak sparsity. Our technique essentially originates from the basic inequality derived from their theoretical analysis of Lasso. The main difference lies that Lasso's results rely on constant lower bounds of some quantities such as the cone invertibility factor or the compatibility factor.   

As for the precision matrix estimation, the error bounds under $\ell_q$ sparsity condition have been discussed for Dantzig-type methods such as the graphical Dantzig method \citep{yuan2010high} and the CLIME method \citep{Cai2011A}, and minimax convergence rates have been established by the ACLIME method \citep{cai2016estimating}. The CLIME method and its variant ACLIME are frequently introduced to deal with the $\ell_q$ sparsity for the precision matrix estimation. Here, our work provides an alternative approach and shows that the Lasso-type method SCIO can obtain the theoretical guarantees of CLIME under the $\ell_q$ sparsity condition. Specially, we relax the irrepresentable condition, which is commonly used for Lasso-type precision matrix estimation. In addition, the SCIO method can be efficiently implemented according to \cite{wang2020efficient} while the computation of Dantzig-type methods turns out to be slow. From this perspective, the SCIO method tends to be more appealing for the high dimensional precision matrix estimation. 

Another closely related Lasso-type method is SLasso proposed by \cite{sun2013sparse}. By inducing a noise level, the SLasso is tuning-free by iteratively estimating the noise level. For the normal distribution, \cite{sun2013sparse} derived the optimal error bounds under the alternative weak sparsity condition, i.e.,  the capped $\ell_1$ measure. The key ideas of SLasso and SCIO are quite similar, e.g.,  the SLasso for fixed noise level $\sigma$ is the same as the SCIO by setting $\hat{\beta}_{jj}=-1$. It would be interesting to compare these two methods from both the computation complexity and the performance of the estimators. We implement the R package ``scalreg" provided by  \cite{sun2013sparse} and it is not very efficient which prevents us from conducting the comparison experiment. From the original paper of \cite{sun2013sparse}, SLasso has improvements over CLIME for most cases. This improvement is due to adaptive choice of the penalty level for each column of the precision matrix. Actually, we can also use different tuning parameters for each column in SCIO (or CLIME) and it is expected to obtain some improvements. Another interesting question is how to exploit the noise level into  complicated cases, e.g., the heavy-tailed data, the non-paranormal data, and the matrix variate data.    We leave these questions as a future work.

For other Lasso-type methods such as the graphical Lasso and D-trace, they are not in a column-by-column form.  Although they have been shown to be consistent under the $\ell_0$ sparsity condition, the extension to the $\ell_q$ sparsity is not trivial and our current technique can not be implemented directly. It is still of interest whether optimal rates can be established under the weak sparse case for the graphical Lasso and D-trace. Moreover, for other statistical problems such as the discriminant analysis problem, the misclassification rate measures the performance of the method and we can use our current technique to derive its error bounds under the general $\ell_q$ sparsity condition. Specifically, it is possible to show that some Lasso-type methods for  the discriminant analysis  such as \cite{fan2012road} or \cite{mai2012direct} are still applicable for the weak sparse case.  We leave these problems for future works.

\section{Acknowledgments}
Wang’s research was supported by the National Natural Science Foundation of China (12031005, 11971017). Liu’s research is supported by National Program on Key Basic Research Project (973 Program, 2018AAA0100704), NSFC Grant No. 12031005, 11825104 and 11690013, Youth Talent Support Program, and a grant from Australian Research Council.

\section{Appendix}
This section includes all the technical proofs of the main theorems and some necessary lemmas.

\subsection{Proof of Theorem \ref{thm:1}}
Let $\be$ be a column of the identity matrix $\mathbf{I}$, then $\bbeta^*=\bSig^{-1} \be$ is the corresponding column of the target precision matrix $\bOme$. For an arbitrary estimator of $\hbSig$, we consider the SCIO estimation
\begin{align*}
	\hbbeta=\underset{\bbeta\in \mR^{p}}{\argmin}\left\{\frac{1}{2} \bbeta\trans \hbSig \bbeta-\be \trans \bbeta+\lambda|\bbeta|_{1}\right\}.
\end{align*}
By the KKT condition, we have
\begin{align} 	
	\hbSig\hbbeta-\be+\lambda\operatorname{sgn}(\hbbeta)=0,
	\label{eq:6}
\end{align}
which ensures the basic inequality
\begin{align}
	(\hbbeta-\bbeta^*)\trans(\hbSig\hbbeta-\be)\leq -\lambda|\hbbeta|_1+\lambda|\bbeta^*|_1.
	\label{eq:7}
\end{align}
Writing the difference vector as $\bh=\bbeta^*-\hbbeta$, we have
\begin{align}
	(\hbbeta-\bbeta^*)\trans(\hbSig\hbbeta-\be)&=(\hbbeta-\bbeta^*)\trans(\hbSig(\hbbeta-\bbeta^*)+(\hbSig-\bSig)\bbeta^*)\notag\\
	&\geq(\hbbeta-\bbeta^*)\trans(\hbSig-\bSig)\bbeta^*\notag\\
	&\geq-|(\hbSig-\bSig)\bbeta^*|_\infty|\bh|_1\notag\\
	&\geq-\|\bOme\|_{L_1}\|\hbSig-\bSig\|_\infty|\bh|_1.
	\label{eq:8}
\end{align}
Combined with the basic inequality \eqref{eq:7}, it reduces to 
\begin{align*}
	-\lambda|\hbbeta|_1+\lambda|\bbeta^*|_1 \geq -\|\bOme\|_{L_1}\|\hbSig-\bSig\|_\infty|\bh|_1.
\end{align*} 	
Since $\lambda\geq 3\|\bOme\|_{L_1}\|\hbSig-\bSig\|_\infty$, then by this assumption we obtain
\begin{align*}
	3\lambda(|\bbeta^*|_1-|\hbbeta|_1)\geq-\lambda|\bh|_1,
\end{align*}
and hence
\begin{align*}
	3(|\bbeta^*|_1-|\hbbeta|_1) \geq-|\bh|_1.
\end{align*}	
For any index set $J$, we have
\begin{align*}
	|\bh_{J^c}|_1& \leq |\bbeta^*_{J^c}|_1+ |\hbbeta_{J^c}|_1= |\bbeta^*_{J^c}|_1+ |\hbbeta|_1-|\hbbeta_{J}|_1\\
	&\leq  |\bbeta^*_{J^c}|_1+|\bbeta^*|_1+\frac{1}{3} |\bh|_1-|\hbbeta_{J}|_1= 2 |\bbeta^*_{J^c}|_1+|\bbeta^*_{J}|_1+\frac{1}{3} |\bh|_1-|\hbbeta_{J}|_1\\
	&\leq 2 |\bbeta^*_{J^c}|_1+|\bh_{J}|_1+\frac{1}{3} |\bh|_1= 2 |\bbeta^*_{J^c}|_1+\frac{4}{3}|\bh_{J}|_1+\frac{1}{3}|\bh_{J^c}|_1,
\end{align*}	
and by rearranging the inequality, we get an important relation 
\begin{align}
	|\bh_{J^c}|_1\leq 2|\bh_J|_1+3|\bbeta^*_{J^c}|_1.
	\label{eq:9}
\end{align}

By the KKT condition \eqref{eq:6}, we have $|\hbSig\hbbeta-\be|_\infty\leq \lambda$ and
\begin{align*}
	|\hbSig{\bbeta^*}-\be|_\infty\leq  \|\hbSig-\bSig\|_\infty|\bbeta^*|_{1} \leq \|\hbSig-\bSig\|_\infty\|\bOme\|_{L_1}\leq\frac{1}{3}\lambda.
\end{align*}
Thus, we can get
\begin{align*}
	|\hbSig\bh|_\infty\leq|\hbSig\hbbeta-\be|_\infty+|\hbSig{\bbeta^*}-\be|_\infty\leq\frac{4}{3}\lambda,
\end{align*}
and
\begin{align*}
	|\bSig\bh|_\infty \leq|(\hbSig-\bSig)\bh|_\infty+|\hbSig\bh|_\infty\leq\|\hbSig-\bSig\|_\infty|\bh|_1+\frac{4}{3}\lambda.
\end{align*}
Then, we conclude that
\begin{align}
	|\bh|_\infty=|\bOme \bSig\bh|_\infty\leq \|\bOme\|_{L_1} 	|\bSig\bh|_\infty & \leq  \|\bOme\|_{L_1}  \|\hbSig-\bSig\|_\infty|\bh|_1+\frac{4}{3} \|\bOme\|_{L_1}  \lambda \nonumber \\
	&\leq  \frac{1}{3} \lambda |\bh|_1+\frac{4}{3} \|\bOme\|_{L_1}  \lambda, 
	\label{eq:10}
\end{align}
and
\begin{align}
	|\bh_J|_1\leq |J| 	|\bh|_\infty \leq  \frac{1}{3} \lambda |J| ( |\bh|_1+4 \|\bOme\|_{L_1}).
	\label{eq:11}
\end{align}

Next, we split the argument into two cases.\\

\noindent \textbf{Case 1:} 	If $|\bh_{J^c}|_1\geq 3|\bh_J|_1$, by the inequality \eqref{eq:9}, we have $|\bh_J|_1\leq 3|\bbeta^*_{J^c}|_1$ and then
\begin{align}
	|\bh|_1=|\bh_{J^c}|_1+ |\bh_{J}|_1 \leq 3|\bh_{J}|_1+3|\bbeta^*_{J^c}|_1\leq 12|\bbeta^*_{J^c}|_1 
	\label{eq:12}
\end{align}
where the first inequality uses the fact \eqref{eq:9} again.

\noindent \textbf{Case 2:} Otherwise, we may assume $|\bh_{J^c}|_1< 3 |\bh_J|_1$ and then $	|\bh|_1 \leq 4|\bh_J|_1$.  By the bound \eqref{eq:11}, 
\begin{align*}
	|\bh|_1 \leq 4|\bh_J|_1 \leq  \frac{4}{3} \lambda |J| ( |\bh|_1+4 \|\bOme\|_{L_1}).
\end{align*}
If the relation $\lambda |J| \leq \frac{1}{2}$ holds, we conclude that
\begin{align}
	|\bh|_1 \leq 16  \lambda |J|   \|\bOme\|_{L_1}.
	\label{eq:13}
\end{align}

Now we begin to conduct the index set $J$ that can control the bounds \eqref{eq:12} and \eqref{eq:13} simultaneously. To do so,  we consider the index set
\begin{align*}
	J=\{j|~|\bbeta_j^*|> t\},
\end{align*}
for some $t>0$. With this setting, 
\begin{align}
	|\bbeta^{*}_{J^c}|_1=\sum_{j \in J^c}|\bbeta^{*}_j|\leq t^{1-q} \sum_{j \in J^c} |\bbeta^{*}_j|^{q}  \leq t^{1-q}s_p,
	\label{eq:14}
\end{align}
and for the cardinality of $J$, we have
\begin{align}
	|J|\leq t^{-q}\sum_{j \in J}\left|\bbeta^{*}_{j}\right|^q \leq t^{-q}s_p.
	\label{eq:15}
\end{align}
Combining the bounds \eqref{eq:12} and \eqref{eq:13}, we get
\begin{align*}
	|\bh|_1  \leq \max\{12|\bbeta^*_{J^c}|_1, 16  \lambda |J|   \|\bOme\|_{L_1} \}\leq \max\{12  t^{1-q}s_p, 16  \lambda t^{-q}s_p   \|\bOme\|_{L_1} \},
\end{align*}
and setting $t=\lambda \|\bOme\|_{L_1}$ yields the conclusion 
\begin{align*}
	|\bbeta^*-\hbbeta|_1\leq16\|\bOme\|_{L_1}^{1-q}\lambda^{1-q}s_p.
\end{align*}
It remains to check $\lambda |J| \leq \frac{1}{2}$. Let $t=\lambda \|\bOme\|_{L_1}$ and we have
\begin{align*}
	\lambda |J|  \leq \lambda  (\lambda \|\bOme\|_{L_1})^{-q} s_p=\lambda^{1-q}  ( \|\bOme\|_{L_1})^{-q} s_p<\frac{1}{2},
\end{align*}	
which holds by the assumption.

Finally, invoking \eqref{eq:10}, we conclude
\begin{align*}
	|\bbeta^*-\hbbeta|_\infty	&\leq\frac{1}{3}\lambda|\bh|_1+\frac{4}{3}\|\bOme\|_{L_1}\lambda\\
	&\leq \frac{16}{3}|\bOme\|_{L_1}^{1-q}\lambda^{2-q}s_p+\frac{4}{3}\|\bOme\|_{L_1}\lambda\\
	&\leq 4\|\bOme\|_{L_1}\lambda,
\end{align*}
where we uses the fact $\lambda^{1-q}  ( \|\bOme\|_{L_1})^{-q} s_p<\frac{1}{2}$ again. 	

\subsection{Proof of Theorem \ref{thm:2}}
Note the result of Theorem \ref{thm:1} holds uniformly for all $i=1,\cdots,p$, that is 
\begin{align*}
	\max_{i=1,\cdots,p}	|\bbeta^*_i-\hbbeta_i|_1\leq16\|\bOme\|_{L_1}^{1-q}\lambda^{1-q}s_p,~\mbox{and}
	\max_{i=1,\cdots,p} 	|\bbeta^*_i-\hbbeta_i|_\infty  \leq 4\|\bOme\|_{L_1}\lambda.
\end{align*}
By the construction of $\tbOme$, it is easy to show
\begin{align*}
	\|\tbOme-\bOme\|_\infty\leq\max_{i=1,\cdots,p} 	|\bbeta^*_i-\hbbeta_i|_\infty  \leq 4\|\bOme\|_{L_1}\lambda.
\end{align*}
Next we study the effect of the symmetrization step. For $i \in \{1,\cdots,n\}$, we denote
\begin{align*}
	\tbbeta=\tbbeta_i,~\hbbeta=\hbbeta_i,\bbeta^*=\bbeta^*_i.
\end{align*}
Since $|\tbbeta|_1 \leq |\hbbeta|_1$ by our construction, for the index set $J$, we have
\begin{align*}
	|\tbbeta_{J^c}|_1=|\tbbeta|_1-|\tbbeta_{J}|_1\leq  |\hbbeta|_1-|\tbbeta_{J}|_1\leq  |(\hbbeta-\tbbeta)_{J}|_1+|\hbbeta_{J^c}|_1,
\end{align*}
which yields 
\begin{align}
	|\tbbeta-\hbbeta|_1\leq |(\hbbeta-\tbbeta)_{J}|_1+|\hbbeta_{J^c}|_1+|\tbbeta_{J^c}|_1\leq 2\{|(\hbbeta-\tbbeta)_{J}|_1+|\hbbeta_{J^c}|_1 \}.
	\label{eq:16}
\end{align}
Recall the index set $J=\{j|~|\bbeta_j^*|> \lambda \|\bOme\|_{L_1}\}$ from Proof of Theorem \ref{thm:1}, and also the bounds
\begin{align*}
	|\bbeta^{*}_{J^c}|_1 \leq ( \lambda \|\bOme\|_{L_1} )^{1-q}s_p,~\mbox{and}~ |J| \leq ( \lambda \|\bOme\|_{L_1})^{-q}s_p.
\end{align*}
Thus
\begin{align*}
	|(\hbbeta-\tbbeta)_{J}|_1\leq |J| |\hbbeta-\tbbeta|_\infty \leq ( \lambda \|\bOme\|_{L_1})^{-q}s_p \cdot 8 \|\bOme\|_{L_1}\lambda=8 ( \lambda \|\bOme\|_{L_1})^{1-q}s_p,
\end{align*}
and 
\begin{align*}
	|\hbbeta_{J^c}|_1\leq |(\hbbeta-\bbeta^*)_{J^c}|_1+|\bbeta^{*}_{J^c}|_1\leq  |\hbbeta-\bbeta^*|_1+|\bbeta^{*}_{J^c}|_1\leq 17 ( \lambda \|\bOme\|_{L_1} )^{1-q}s_p.
\end{align*}	
Invoking the bound \eqref{eq:16}, we get
\begin{align*}
	|\tbbeta-\hbbeta|_1\leq 50 ( \lambda \|\bOme\|_{L_1} )^{1-q}s_p,
\end{align*}	
which ensures 
\begin{align*}
	|\tbbeta-\bbeta^*|_1\leq |\tbbeta-\hbbeta|_1+|\hbbeta-\bbeta^*|_1\leq 66 ( \lambda \|\bOme\|_{L_1} )^{1-q}s_p.
\end{align*}
Since the above bound holds uniformly for all $i=1,\cdots,p$, we conclude
\begin{align*}
	\|\tbOme-\bOme\|_{L_1} \leq \max_{i=1,\cdots,p}|\tbbeta-\bbeta^*|_1\leq 66 ( \lambda \|\bOme\|_{L_1} )^{1-q}s_p.
\end{align*}	 

\subsection{Proof of  Corollaries}
We only prove Corollary 1. The proof of other corollaries share a very similar procedure as Corollary 1 and hence are omitted.
\bigskip

\noindent \textbf{Proof of Corollary 1:}\\
By the assumption (A2) and Proposition \ref{prop:1}, we have $\mP(\|\hbSig_1-\bSig\|_{\infty}\geq C\sqrt{\frac{\log{p}}{n}})= O(p^{-\tau})$. Similarly, by the assumption (A3) and Proposition \ref{prop:1}, we have $\mP(\|\hbSig_1-\bSig\|_{\infty}\geq C\sqrt{\frac{\log{p}}{n}})= O(p^{-\tau}+n^{-\frac{\delta}{8}})$.

We take $\lambda=C_0M_p\sqrt{\frac{\log{p}}{n}}$. Note that $\|\bOme\|_{L_1}\leq M_p$, then for a sufficiently large $C_0$, the condition $\lambda\geq3\|\bOme\|_{L_1}\|\hbSig-\bSig\|_\infty$ required in Theorem \ref{thm:2} holds.
By the assumption (A1), when $n,p$ are large enough, we can obtain that $\|\bOme\|^{-q}_{L_1}\lambda^{1-q}s_p\leq s_pM_p^{1-2q}\leq \frac{1}{2}$. So by applying Theorem \ref{thm:2}, the conclusion of Corollary \ref{cor:1} holds.

\bigskip

\bibliographystyle{chicago}      
\bibliography{mybibfile}

\end{document}